# Homogeneity of infinite dimensional isoparametric submanifolds

By Ernst Heintze and Xiaobo Liu

A subset $S$ of a Riemannian manifold $N$ is called *extrinsically homogeneous* if $S$ is an orbit of a subgroup of the isometry group of $N$. In [Th], Thorbergsson proved the remarkable result that *every complete, connected, full, irreducible isoparametric submanifold of a finite dimensional Euclidean space of rank at least 3 is extrinsically homogeneous.* This result, combined with results of [PT1] and [Da], finally classified irreducible isoparametric submanifolds of a finite dimensional Euclidean space of rank at least 3. While Thorbergsson's proof used Tits buildings, a simpler proof without using Tits buildings was given by Olmos (cf. [O2]). The main purpose of this paper is to extend Thorbergsson's result to the infinite dimensional case.

The study of infinite dimensional isoparametric submanifolds of a Hilbert space (always assumed to be separable) was initiated by Terng [T2]. Besides its intrinsic interest, the theory of infinite dimensional isoparametric submanifolds is a very useful tool in studying the submanifold geometry of compact symmetric spaces (cf. [TT] as well as [HL] and [Ew]). Although much progress has been made (especially, many basic properties of finite dimensional isoparametric submanifolds having been successfully extended to the infinite dimensional case (cf. [T2], [T5], and [HL])), the classification of infinite dimensional isoparametric submanifolds is far from being solved. To this end, the understanding of the homogeneity of infinite dimensional isoparametric submanifolds will certainly play a very important role, as it does in the finite dimensional case. Examples of non-homogeneous infinite dimensional isoparametric hypersurfaces have been found by Terng and Thorbergsson (cf. [TT]). In this paper, we will prove the following theorem which solves a long standing open problem (cf. [T3], [T4], [TT]).

THEOREM A. *Let $M$ be a complete, connected, irreducible isoparametric submanifold in a Hilbert space $V$. Assume that the set of all the curvature normals of $M$ at some point is not contained in any affine line. Then $M$ is extrinsically homogeneous in the Hilbert space $V$.*

*Remark.* In this theorem, the dimension of $M$ and thus $V$ could be either finite or infinite. Without loss of generality, we may assume that $M$ is full, i.e.

not contained in any proper affine subspace. Then the assumption that all the curvature normals of $M$ at some point do not lie in any affine line is equivalent to the condition that the codimension of $M$ is bigger than or equal to 2 if $M$ is infinite dimensional or the codimension of $M$ is bigger than or equal to 3 if $M$ is finite dimensional (cf. [T1], [T2] as well as [HL]). Therefore our proof, when applied to the finite dimensional case, also gives a new approach to Thorbergsson's theorem which simplifies the arguments of Thorbergsson and Olmos. Our proof differs from Thorbergsson's proof in that we do not use Tits buildings. It differs from Olmos' proof in that we construct the extrinsic isometries without constructing a canonical connection. Note that in the finite dimensional case, the curvature normals at one point are always contained in an affine hyperplane of the normal space of $M$ at that point, whereas in infinite dimensions this is not the case.

In the proof of the above theorem, we also obtained some other results which are interesting in their own right. To state those results more precisely, we introduce some notation first. Recall that for a submanifold $M$ of a Hilbert space $V$, the *end point map* $\eta : \nu M \longrightarrow V$ is defined by $\eta(v) = x + v$ for $v \in \nu M_x$, where $\nu M$ is the normal bundle of $M$. $M$ is called *proper Fredholm* if it has finite codimension and the end point map restricted to any finite normal disk bundle is a proper Fredholm map. A proper Fredholm submanifold $M$ is called *isoparametric* if its normal bundle is globally flat, and the shape operators along any parallel normal vector field are conjugate. We will always assume that $M$ is complete and connected. Since the normal bundle $\nu M$ is flat, the tangent bundle of $M$ splits as $TM = \overline{\bigoplus \{E_i \mid i \in I\}}$ into the direct sum of the simultaneous eigenspaces $E_i$ of the shape operators, where $I$ is a countable index set. In the set $\{E_i \mid i \in I\}$, the $E_i$ are called the curvature distributions of $M$. Let $\{n_i \mid i \in I\}$ be the corresponding curvature normals of $M$, i.e. the globally defined parallel normal vector fields such that for any parallel normal vector field $v$ on $M$, the restriction of the shape operator to each $E_i$ is

$$A_v|_{E_i} = \langle v, n_i \rangle \operatorname{Id}.$$

We will always denote the zero curvature normal by $n_0$ (if it occurs) and the corresponding curvature distribution by $E_0$. Furthermore we will always assume that $M$ is full. This is equivalent to saying that for any point $x \in M$ the curvature normals $\{n_i(x) \mid i \in I\}$ span the normal space $\nu_x M$. It is known that each curvature distribution is integrable. If $n_i \neq 0$, the rank of $E_i$ is finite and for any $x \in M$, the leaf of $E_i$ passing through $x$, denoted by $S_i(x)$, is a round sphere centered at $x + (n_i(x)/\|n_i\|^2)$ with radius $1/\|n_i\|$. Also $S_i(x)$ is called the $i^{\text{th}}$ curvature sphere at $x$. The leaves of $E_0$ are closed affine subspaces of the Hilbert space $V$ (cf. [T2]). Two points $p$ and $q$ in $M$ are called *equivalent*, denoted by $p \sim q$, if there exists a sequence of points $x_j$, $j = 0, \ldots, n$, such that $x_0 = p$, $x_n = q$, and for every $j = 1, \ldots, n$, $x_{j-1}$ and $x_j$ lie in one curvature sphere. For any $p \in M$, the equivalence class of $p$, denoted by $Q(p)$, is defined to be
$$Q(p) = \{q \in M \mid q \sim p\}.$$
One key observation is that each equivalence class contains almost all the information of $M$. In fact, if $M$ is a compact, connected isoparametric submanifold of a finite dimensional Euclidean space, then each equivalence class is equal to $M$ (cf. [HOTh, Lemma 3.1]). In the infinite dimensional case, we have the following result which is nontrivial even for the known homogeneous examples.

THEOREM B. *Let $M$ be an irreducible isoparametric submanifold in a Hilbert space $V$ with codimension at least 2. Then for any $p \in M$, $Q(p)$ is dense in $M$.*

This result is a crucial step in the proof of Theorem A. Another crucial step is that, under the same assumption as in Theorem B, each equivalence class of $M$ is actually extrinsically homogeneous.

One immediate consequence of Theorem A is the homogeneous slice theorem for infinite dimensional isoparametric submanifolds. Fix one point $x_0 \in M$. Given an affine subspace $P$ of the normal space $\nu_{x_0} M$, there is a corresponding distribution on $M$ defined by
$$D_P = \overline{\bigoplus \{E_i \mid n_i(x_0) \in P\}}.$$
If $P' \subset \nu_{x_0} M$ is another affine subspace which contains the same curvature normals as $P$ does, then $D_P = D_{P'}$. Thus there is a bijection between the affine subspaces $P$ of $\nu_{x_0} M$ which are spanned by curvature normals and the distribution $D_P$. Note that in case $M$ is compact and hence finite dimensional, all curvature normals are contained in an affine subspace which does not contain 0. It is known that the distribution $D_P$ is always integrable and has totally geodesic leaves (cf. [HL, Prop. 2.3]). Let $L_P(x)$ be the leaf of $D_P$ passing through $x \in M$, and define $W_P(x) = x + D_P(x) \oplus \mathrm{span}\{n_i(x) \mid n_i(x_0) \in P\}$. Then $L_P(x)$ is a full isoparametric submanifold of $W_P(x)$, and $W_P(y) = W_P(x)$ for all $y \in L_P(x)$ (cf. [HL, Lemma 3.3]). We call $L_P(x)$ the *slice* at $p$ corresponding to $D_P$. The following result is an extension of the homogeneous slice theorem (cf. [HOTh]) to the infinite dimensional case.

THEOREM C. *If $M$ is a full, irreducible isoparametric submanifold of an infinite dimensional Hilbert space with codimension at least 2, then $L_P(x)$ is extrinsically homogeneous in $W_P(x)$ for any affine subspace $P$ of $\nu_{x_0} M$ and any $x \in M$.*

Notice that in the finite dimensional case, the homogeneous slice theorem has been used in both proofs of Thorbergsson's theorem by Thorbergsson and Olmos. However, in the infinite dimensional case, the homogeneous slice theorem is more difficult to prove and is one of the major obstacles to generalize Olmos' proof of Thorbergsson's theorem to infinite dimensions. In this paper, we will use a weak version of Theorem C, i.e. the homogeneity of finite dimensional slices, to prove Theorem A. Theorem C will be proved as a consequence of Theorem A.

In the proofs of the theorems above, we need a generalization of a theorem of Chow (cf. [Ch]) to infinite dimensions. Let $N$ be a complete connected Hilbert manifold. Let $D$ be a set of smooth vector fields which are defined on open subsets of $N$. Two points $x$ and $y$ are called $D$-equivalent, denoted by $x \sim_D y$, if $x$ and $y$ can be connected by a piecewise differentiable curve, each differentiable piece of which is an integral curve of a vector field in $D$. Let $\Omega_D(x) = \{y \in N \mid y \sim_D x\}$. We call $\Omega_D(x)$ the set of reachable points of $D$ starting from $x$. Let $D^*$ be the set of smooth vector fields on $N$ which is generated by $D$ in the following sense: $D \subset D^*$, $D^*$ contains the zero vector field, and for any $X$ and $Y$ belonging to $D^*$, so do the vector fields $aX + bY$ and $[X, Y]$, which are defined on the intersection of the domains of $X$ and $Y$, for any real numbers $a$ and $b$. For any $x \in N$, let $D^*(x)$ be the set of all vectors $X(x)$ where $X$ is a vector field defined in a neighborhood of $x$ and $X \in D^*$. If $N$ is finite dimensional and $D^*(x) = T_x N$ for all $x \in N$, then Chow's theorem says that $\Omega_D(x) = N$ for all $x \in N$. A generalization of Chow's theorem for arbitrary sets of vector fields was given by Sussmann ([Su]) in the finite dimensional case. It is known that those theorems have important applications in control theory. In infinite dimensions, the following version of Chow's theorem is true.

THEOREM D (generalized Chow's theorem). *Let $N$ be a complete connected Hilbert manifold and $D$ a set of smooth vector fields defined on open subsets of $N$. If $D^*(x)$ is dense in $T_x N$ for all $x \in N$, then $\Omega_D(x)$ is dense in $N$ for all $x \in N$.*

We do not know whether this result has been proved before. Since we cannot find relevant references, we include a proof in the appendix of this paper for completeness.

We will mainly focus on the proof of Theorem A for the infinite dimensional case. The proof for finite dimensions can be done similarly. Many arguments needed for the infinite dimensional case are redundant for the finite dimensional case, in particular the results of Sections 1, 2, 5, and the appendix. This paper is organized as follows. In Section 1, we define the notion of horizontal curves and prove that associated to each horizontal curve, there is an isometry between the corresponding slices (which could be either finite or infinite dimensional). Such isometries play a crucial role in the construction of the global isometries of equivalence classes. In Section 2, we prove the homogeneous slice theorem for finite dimensional slices and some properties of the normal holonomy groups of focal submanifolds. In Section 3, we give a criterion for two isoparametric submanifolds to coincide with each other. This criterion and the properties for the normal holonomy groups will be used in Section 4 to prove that the isometry we constructed preserves the isoparametric submanifold. In Section 4, we construct extrinsic isometries of $M$ and prove the homogeneity of the closure of an equivalence class, which implies Thorbergsson's theorem for the finite dimensional case. Theorems A to C are proved in Section 5. Finally we give a proof of Theorem D in the appendix.

We would like to thank Jens Heber for showing us a proof of Lemma A.1 in the appendix. The second author would like to thank the University of Augsburg and both authors would like to thank the Max-Planck-Institut für Mathematik at Bonn for their support and hospitality.

# 1. Horizontal curves

In this section, we will define the notion of horizontal curves and prove some of their properties. This notion is motivated by [HOTh] (see also the remark below). Given an integrable distribution $D$ on $M$, we call a piecewise differentiable curve $\gamma$ in $M$ *horizontal* with respect to $D$ if $\dot\gamma(t) \perp D(\gamma(t))$ for all $t$ where $\gamma$ is differentiable. A horizontal curve $\gamma$ is called *parallel* to a curve $\beta$ in $M$ with respect to $D$ if for every $t$, $\gamma(t)$ and $\beta(t)$ always lie in the same leaf of $D$.

Fix one point $x_0 \in M$. Let $P$ be an affine subspace of the normal space $\nu_{x_0}M$, and $D_P$ the corresponding distribution on $M$ defined in the introduction. The main purpose of this section is to prove the following

PROPOSITION 1.1. *Let $\gamma$ be a horizontal curve with respect to $D_P$. Suppose that the domain of $\gamma$ is an interval which contains $0$. Then there exists a one-parameter family of isometries $f_t : W_P(\gamma(0)) \longrightarrow W_P(\gamma(t))$, defined for all $t$ in the domain of $\gamma$, such that $f_t(L_P(\gamma(0))) = L_P(\gamma(t))$ and for any $x \in L_P(\gamma(0))$, the curve $\gamma_x(t) := f_t(x)$ is the unique horizontal curve through $x$ which is parallel to $\gamma$.*

*Remarks.* (i) In case that $P$ is an affine subspace which does not contain $0$, Proposition 1.1 follows easily, even if $\gamma$ is not horizontal. In fact, there exists a parallel normal vector field $\xi$ on $M$ such that for any curvature normal $n_i$, $\langle n_i(x_0), \xi(x_0)\rangle = 1$ if and only if $n_i(x_0) \in P$. Let $\eta$ be the end-point map, i.e. $\eta(v) = x + v$ for all $x \in M$ and $v \in \nu_x M$. The set $M_\xi := \{\eta(\xi(x)) \mid x \in M\}$ is an embedded proper Fredholm submanifold of $V$ and $\pi_\xi : M \longrightarrow M_\xi$, which is defined by $\pi_\xi(x) = \eta(\xi(x))$, is a fibration with finite dimensional fibers. Furthermore $\ker(\pi_\xi)_*|_x = D_P(x)$. The connected com-ponent of $(\pi_\xi)^{-1}\pi_\xi(x)$ passing through $x$ is $L_P(x)$. Let $\gamma_\xi(t) = \pi_\xi(\gamma(t))$. Then $\gamma_\xi$ is a piecewise differentiable curve in $M_\xi$ and $W_P(\gamma(t)) \subset \gamma_\xi(t) + \nu_{\gamma_\xi(t)} M_\xi$. Let $\tau_t : \nu_{\gamma_\xi(0)} M_\xi \longrightarrow \nu_{\gamma_\xi(t)} M_\xi$ be the parallel translation along $\gamma_\xi$ in the normal bundle of $M_\xi$. Then $f_t(x) := \gamma_\xi(t) + \tau_t(x - \gamma_\xi(0))$ apparently defines an isometry from $W_P(\gamma(0))$ to $W_P(\gamma(t))$ which satisfies the properties in the above proposition.

(ii) Assume in addition that $\gamma$ is defined over the unit interval $[0, 1]$ and that $\pi_\xi(\gamma(1)) = \pi_\xi(\gamma(0))$. Then $f_1$ is an isometry from $W_P(\gamma(0))$ to itself. The group generated by all such isometries is exactly the normal holonomy group of $M_\xi$ at $\gamma_\xi(0)$.

Before proving the proposition, we first notice the following general fact.

LEMMA 1.2. *Let $\alpha_s$ be a one-parameter family of differentiable curves in a Riemannian manifold $N$ which satisfies the property that for every $t$, there exists a totally geodesic submanifold $T_t$ of $N$ such that $\alpha_s(t) \in T_t$ and $\frac{\partial}{\partial t}\alpha_s(t) \perp T_t$ for all $s$. Assume that $\alpha$ is differentiable in $(s, t)$. Then for every $s$, $\left\|\frac{\partial}{\partial s}\alpha_s(t)\right\|$ does not depend on $t$.*

*Proof.*
$$\frac{\partial}{\partial t}\left\|\frac{\partial}{\partial s}\alpha_s(t)\right\|^2 = 2\left\langle \frac{D}{\partial t}\frac{\partial}{\partial s}\alpha_s(t), \frac{\partial}{\partial s}\alpha_s(t)\right\rangle = 2\left\langle \frac{D}{\partial s}\frac{\partial}{\partial t}\alpha_s(t), \frac{\partial}{\partial s}\alpha_s(t)\right\rangle.$$

By assumption, $\frac{\partial}{\partial s}\alpha_s(t)$ is always tangent to $T_t$ and $\frac{\partial}{\partial t}\alpha_s(t)$ is always perpendicular to $T_t$. Therefore
$$\frac{\partial}{\partial t}\left\|\frac{\partial}{\partial s}\alpha_s(t)\right\|^2 = -2\left\langle \frac{\partial}{\partial t}\alpha_s(t), \frac{D}{\partial s}\frac{\partial}{\partial s}\alpha_s(t)\right\rangle.$$

Since $T_t$ is totally geodesic, $\frac{D}{\partial s}\frac{\partial}{\partial s}\alpha_s(t)$ is always tangent to $T_t$. Hence $\frac{\partial}{\partial t}\left\|\frac{\partial}{\partial s}\alpha_s(t)\right\|^2 \equiv 0$. This proves the lemma. □

*Proof of Proposition* 1.1. Based on the remark following the proposition, we only need to consider the case where $P$ is a linear subspace of $\nu_{x_0} M$. Notice that in this case, $D_P$ contains $E_0$ as a subdistribution since we always treat $n_0 = 0$ as a curvature normal. We may assume that $P$ is spanned by curvature normals which belong to $P$ because otherwise we can work with a proper linear subspace of $P$ which satisfies this condition. Without loss of generality, we may also assume that $\gamma$ is a regular differentiable horizontal curve with respect to $D_P$ and is defined over the unit interval $[0, 1]$. The general case can be handled by breaking $\gamma$ into a finite number of regular differentiable pieces and using induction on their number. Choose a parallel normal vector field $\xi$ on $M$ such that for every curvature normal $n_i$, $n_i(x_0) \perp \xi(x_0)$ if and only if $n_i(x_0) \in P$. Let $F : M \longrightarrow V$ be the Gauss map defined by $\xi$, i.e., $F(x) = \xi(x)$. Since $F_* |_x = -A_{\xi(x)}$, where $x \in M$ and $A$ is the shape operator of $M$, we have that $\ker(F_* |_x) = D_P(x)$ and the connected component of $F^{-1}(F(x))$ containing $x$ is $L_P(x)$. Denote the image of $\gamma$ under the map $F$ by $\bar{\gamma}$. Then $\bar{\gamma}$ is a regular curve since $\gamma$ is regular and horizontal. Without loss of generality, we may further assume that $\bar{\gamma}$ does not have self-intersections. The proof is divided into three steps.

*Step* 1. For any $x \in M$ and $r > 0$, let $B_r(x)$ be the open distance ball in $L_P(x)$ with radius $r$ and center $x$. There exists an $\varepsilon > 0$ such that through every $x \in B_\varepsilon(\gamma(0))$, there exists a unique regular differentiable horizontal curve $\gamma_x$, defined over $[0, 1]$, which is parallel to $\gamma$.

We first prove that there exists $\varepsilon' > 0$ such that the set $N_{\varepsilon'} := \bigcup_{t \in [0,1]} B_{\varepsilon'}(\gamma(t))$ is an immersed submanifold of $M$ with boundary. In fact, for every $t_0 \in [0, 1]$, there is a neighborhood of $\gamma(t_0)$ which is diffeomorphic to $V_1 \times V_2$, where $V_1$ and $V_2$ are open subsets of two Hilbert spaces, and the leaves of $D_P$, when restricted to this neighborhood, are given by $V_1 \times \{y_2\}$, where $y_2 \in V_2$. The curve $\gamma$, when restricted to this neighborhood, can be written as $\gamma(t) = (\gamma_1(t), \gamma_2(t))$, where $\gamma_i(t) \in V_i$ for $i = 1, 2$. Since $\gamma$ is regular and horizontal with respect to $D_P$, the curve $\gamma_2$ is a regular differentiable curve in $V_2$. Therefore the set $V_1 \times \gamma_2$ is an immersed submanifold of $V_1 \times V_2$. It follows that there is a neighborhood of $t_0$ and a positive number $\varepsilon_{t_0}$ such that $\bigcup B_{\varepsilon_{t_0}}(\gamma(t))$, where $t$ runs over this neighborhood of $t_0$, is an immersed submanifold of $M$ (with boundary if $t_0 = 0$ or $1$). Since $[0, 1]$ is compact, we can choose $\varepsilon' > 0$ such that $N_{\varepsilon'}$ is an immersed submanifold of $M$ with boundary.

When restricted to $N_{\varepsilon'}$, $D_P$ is a codimension one distribution. Therefore, we can define a unit vector field $u$ on $N_{\varepsilon'}$ such that for every $x \in N_{\varepsilon'}$, $u(x)$ is tangent to $N_{\varepsilon'}$, $u(x) \perp D_P(x)$, and $\langle F_*(u(x)), \dot{\bar{\gamma}} |_{F(x)} \rangle > 0$. Since the image of $N_{\varepsilon'}$ under the Gauss map $F$ is $\bar{\gamma}$, $F_*(u(x))$ is always proportional to $\dot{\bar{\gamma}} |_{F(x)}$. Moreover, $F_*(u(x)) \neq 0$ since $\ker(F_* |_x) = D_P$. Therefore $u$ is well defined. For any $x \in N_{\varepsilon'}$, let $g(x) = \|F_*(u(x))\|$. Since $g(x) > 0$ and $g$ is smooth on $N_{\varepsilon'}$, by the compactness of the unit interval, there exist positive numbers $a$ and $\varepsilon < \varepsilon'$ such that $g(x) > a$ for all $x \in N_\varepsilon := \bigcup_{t \in [0,1]} B_\varepsilon(\gamma(t))$. For any $x \in N_\varepsilon$, let $\pi(x)$ be the unique point on $\gamma$ such that $x$ and $\pi(x)$ belong to the same connected component of $F^{-1}(F(x))$. Now $\pi(x)$ is well defined since $\bar{\gamma}$ has no self-intersections. Let $w$ be the vector field on $N_\varepsilon$ which is defined by $w(x) = \{\|\dot{\gamma} |_{\pi(x)}\| \cdot g(\pi(x))/g(x)\} \cdot u(x)$. Then $w$ is a bounded vector field on $N_\varepsilon$ such that $F_*(w(x)) = \dot{\bar{\gamma}} |_{F(x)}$ for all $x$ and $w(\gamma(t)) = \dot{\gamma}(t)$ for all $t$. For any $x \in B_\varepsilon(\gamma(0))$, let $\gamma_x$ be the integral curve of $w$ which passes through $x$. It follows from the definition of $w$ that, on the interval where it is defined, $\gamma_x$ is a horizontal curve parallel to $\gamma$. It remains to show that $\gamma_x$ is defined over the whole interval $[0, 1]$. Note that the uniqueness of the horizontal curve which is parallel to $\gamma$ and passing through $x$ follows from the uniqueness of the integral curve of $w$ passing through $x$.

For any $x \in B_\varepsilon(\gamma(0))$, there exists a piecewise differentiable curve $\beta$ in $B_\varepsilon(\gamma(0))$ such that $\beta(0) = \gamma(0)$, $\beta(1) = x$ and the length of $\beta$ is less than $\varepsilon$. Let $t_0$ be the largest positive number such that $\gamma_{\beta(s)}$ is defined over $[0, t_0)$ for all $s \in [0, 1]$. Then $t_0$ exists since the image of $\beta$ is compact. For any $s$, since $w$ is bounded, $\gamma_{\beta(s)}$ converges (uniformly in $s$) to a point in $M$ which is denoted by $\gamma_{\beta(s)}(t_0)$. We have $\gamma_{\beta(s)}(t_0) \in L_P(\gamma(t_0))$ since $F(\gamma_{\beta(s)}(t_0)) = \bar{\gamma}(t_0)$ and $s \longmapsto \gamma_{\beta(s)}(t_0)$ is a continuous curve starting at $\gamma(t_0)$. By Lemma 1.2, for each $t \in [0, t_0)$, the length of the curve $s \longmapsto \gamma_{\beta(s)}(t)$ is equal to the length of $\beta$, which is less than $\varepsilon$. Therefore $\gamma_{\beta(s)}(t_0) \in B_\varepsilon(\gamma(t_0))$. So for all $s$, $\gamma_{\beta(s)}$ is defined over $[0, t_0]$. If $t_0 < 1$, by the compactness of the curve $s \longmapsto \gamma_{\beta(s)}(t_0)$, all the curves $\gamma_{\beta(s)}$ can be extended further to an open interval which contains $[0, t_0]$. This contradicts the definition of $t_0$. Therefore $t_0 = 1$. This proves the desired statement.

*Step* 2. There is a one-parameter family of isometries $f_t : W_P(\gamma(0)) \longrightarrow W_P(\gamma(t))$, $t \in [0, 1]$, such that $f_t(B_\varepsilon(\gamma(0))) \subset B_\varepsilon(\gamma(t))$ and for any $x \in B_\varepsilon(\gamma(0))$, the curve $\gamma_x(t) := f_t(x)$ is the unique horizontal curve through $x$ which is parallel to $\gamma$.

Recall that there is a canonical way to associate a Coxeter group to $M$ which acts on the affine normal space $x_0 + \nu_{x_0} M$. Each connected component of the set of regular points of this action is called an open Weyl chamber (cf. [T2]). Let $\Delta$ be an open Weyl chamber in $x_0 + \nu_{x_0} M$ containing $x_0$. For $x \in M$, let $P(x)$ and $\Delta(x)$ be the subsets of $\nu_x M$ which are obtained by parallel translating $P$ and $\Delta$, respectively, from $x_0$ to $x$ in the normal bundle. For $t \in [0, 1]$, let $U_\varepsilon(t) = \{x + v \mid x \in B_\varepsilon(\gamma(t)),\ v \in P(\gamma(t)) \cap \Delta(\gamma(t))\}$. Then $U_\varepsilon(t)$ is an open subset of $W_P(\gamma(t))$. For any $y = x + v \in U_\varepsilon(0)$, where $x \in B_\varepsilon(\gamma(0))$ and $v \in P(\gamma(0)) \cap \Delta(\gamma(0))$, let $\gamma_y(t) = \gamma_x(t) + v(\gamma_x(t))$, where $\gamma_x$ is defined by Step 1 and $v(\gamma_x(t))$ is the parallel translation of $v$ from $x$ to $\gamma_x(t)$ in the normal bundle. Since $\gamma_x(t) \in B_\varepsilon(\gamma(t))$, $\gamma_y(t) \in U_\varepsilon(t) \subset W_P(\gamma(t))$. Moreover $\dot{\gamma}_y(t) = \dot{\gamma}_x(t) - A_v(\dot{\gamma}_x(t)) \perp W_P(\gamma(t))$. Since for each $t$, $W_P(\gamma(t))$ is a totally geodesic submanifold of the Hilbert space $V$, by Lemma 1.2, the map $f_t : U_\varepsilon(0) \longrightarrow U_\varepsilon(t)$ defined by $f_t(y) = \gamma_y(t)$ is an isometry. Since $U_\varepsilon(t)$ is a connected open subset of the affine subspace $W_P(\gamma(t))$, this isometry can be extended in a unique way to an isometry from $W_P(\gamma(0))$ to $W_P(\gamma(t))$. By abuse of notation, we still denote it by $f_t$. Then $f_t$ apparently satisfies the required conditions.

*Step* 3. $f_t(L_P(\gamma(0))) = L_P(\gamma(t))$ and for any $x \in L_P(\gamma(0))$, the curve $\gamma_x(t) := f_t(x)$ is the unique horizontal curve through $x$ which is parallel to $\gamma$.

If for some $x \in L_P(\gamma(0))$, $\gamma_x$ is a horizontal curve in $M$ which is parallel to $\gamma$, then by the first two steps, there is a one-parameter family of isometries, denoted by $g_t^x$, from $W_P(\gamma(0))$ to $W_P(\gamma(t))$ and an open neighborhood of $x$, denoted by $B(x)$, in $L_P(\gamma(0))$ such that $g_t^x(B(x)) \subset L_P(\gamma(t))$ and for any $z \in B(x)$, the curve $t \longmapsto g_t^x(z)$ is the unique horizontal curve through $z$ which is parallel to $\gamma$. We call $g_t^x$ the one-parameter family of isometries defined by $\gamma_x$. For any $x \in L_P(\gamma(0))$, we choose an arbitrary curve $\beta$ in $L_P(\gamma(0))$ such that $\beta(0) = \gamma(0)$ and $\beta(1) = x$. Let $a$ be the largest positive number satisfying the conditions that for all $s \in [0, a)$, $\gamma_{\beta(s)}$ is a horizontal curve in $M$ which is parallel to $\gamma$ and the one-parameter family of isometries $g_t^{\beta(s)}$ defined by $\gamma_x$ coincides with $f_t$. Such number exists by Step 2. Since $f_t$ is an isometry between two closed affine subspaces, $\gamma_{\beta(s)}(t)$ converges to $\gamma_{\beta(a)}(t)$ and $\dot{\gamma}_{\beta(s)}(t)$ converges to $\dot{\gamma}_{\beta(a)}(t)$ as $s$ goes to $a$ from below. Since $M$ is complete, $\gamma_{\beta(a)}$ is contained in $M$ and is horizontal with respect to $D_P$. Since $F(\gamma_{\beta(a)}(t)) = F(\gamma(t))$, where $F$ is the Gauss map, and $\gamma_{\beta(a)}(t)$ can be connected to $\gamma(t)$ by the curve $s \longmapsto \gamma_{\beta(s)}(t)$, $\gamma_{\beta(a)}(t) \in L_P(\gamma(t))$. Therefore $\gamma_{\beta(a)}$ is parallel to $\gamma$. To prove that the one-parameter family of isometries $g_t^{\beta(a)}$ defined by $\gamma_{\beta(a)}$ coincides with $f_t$, we choose an arbitrary $s \in [0, a)$ such that $\beta(s) \in B(\beta(a))$. For any $z \in B(\beta(s)) \cap B(\beta(a))$, both curves $t \longmapsto f_t(z)$ and $t \longmapsto g_t^{\beta(a)}(z)$ are horizontal curves through $z$ which are parallel to $\gamma_{\beta(a)}$. By the uniqueness of parallel horizontal curves, we have $f_t(z) = g_t^{\beta(a)}(z)$. This implies that $f_t$ has to agree with $g_t^{\beta(a)}$ over the whole $W_P(\gamma(0))$ since $L_P(\gamma(0))$ is full in $W_P(\gamma(0))$. If $a < 1$, we can choose $b > a$ such that $\beta([a, b]) \subset B(\beta(a))$. The number $b$ satisfies the conditions in the definition of $a$. This contradicts the choice of $a$. Therefore $a = 1$. It follows that $\gamma_x$ is a horizontal curve parallel to $\gamma$ for all $x \in L_P(\gamma(0))$ and $f_t(L_P(\gamma(0))) \subset L_P(\gamma(t))$. The fact that $f_t(L_P(\gamma(0))) \supset L_P(\gamma(t))$ follows by consideration of the one-parameter family of isometries induced by the curve $t \longmapsto \gamma(1-t)$ which gives the inverse of $f_t$. The proof of the proposition is thus finished. $\square$

In the rest of this section, we prove some properties of the one parameter family of isometries defined by a horizontal curve.

LEMMA 1.3. *Let $\gamma$ be a horizontal curve with respect to $D_P$. Let $f_t$ be the one-parameter family of isometries defined in Proposition 1.1. For any $x \in L_P(\gamma(0))$, $(f_t)_* |_{\nu_x L_P(\gamma(0))}$ is the restriction of the parallel translation in $\nu M$ to $\nu_x L_P(\gamma(0))$, where $\nu_x L_P(\gamma(0))$ should be understood as a subspace of $W_P(\gamma(0))$. Moreover for any curvature normal $n_i(x_0) \in P$, $(f_t)_*(E_i(x)) = E_i(f_t(x))$.*

*Proof.* The first assertion follows from the construction of $f_t$ in Step 2 of the proof of Proposition 1.1. Since $f_t$ is an isometry from $W_P(\gamma(0))$ to $W_P(\gamma(t))$, it maps the curvature sphere $S_i(x) \subset L_P(\gamma(0))$ to some curvature sphere of $L_P(\gamma(t))$ with curvature normal $(f_t)_*(n_i(x))$. By the first assertion of this lemma, $(f_t)_*(n_i(x)) = n_i(f_t(x))$. Therefore $f_t(S_i(x)) = S_i(f_t(x))$. This proves the second assertion of the lemma. □

If $P_1$ and $P_2$ are two affine subspaces of $\nu_{x_0} M$ such that $P_1 \subset P_2$, and $\gamma$ is a horizontal curve with respect to $D_{P_2}$, then $\gamma$ is also horizontal with respect to $D_{P_1}$, which is a subdistribution of $D_{P_2}$. Let $f_t^i : W_{P_i}(\gamma(0)) \longrightarrow W_{P_i}(\gamma(t))$, $i = 1, 2$, be the corresponding one-parameter family of isometries defined by Proposition 1.1.

LEMMA 1.4.  $f_t^2(L_{P_1}(\gamma(0))) = L_{P_1}(\gamma(t))$ and $f_t^2 |_{W_{P_1}(\gamma(0))} = f_t^1$.

*Proof.* By Lemma 1.3, $(f_t^2)_*$ preserves the distribution $D_{P_1}$. Therefore $f^2$ always maps one leaf of $D_{P_1}$ to another leaf of $D_{P_1}$. This proves the first equality. Now for any $x \in L_{P_2}(\gamma(0))$, let $\gamma_x$ be the horizontal curve with respect to $D_{P_2}$ which starts at $x$ and is parallel to $\gamma$. Since $D_{P_1} \subset D_{P_2}$, $\gamma_x$ is also horizontal with respect to $D_{P_1}$. Moreover, if $x \in L_{P_1}(\gamma(0))$, then $\gamma_x(t) \in L_{P_1}(\gamma(t))$ by the first equality of this lemma. Therefore $\gamma_x$ is also parallel to $\gamma$ as a horizontal curve with respect to $D_{P_1}$. By Proposition 1.1, $f_t^2(x) = \gamma_x(t) = f_t^1(x)$. Since $L_{P_1}(\gamma(0))$ is full in $W_{P_1}(\gamma(0))$, this proves the second equality. □

## 2. Normal holonomy groups

In [O1], Olmos proved that the restricted normal holonomy group of a submanifold of a finite dimensional space of constant curvature is compact and the nontrivial part of its representation on the normal space is an s-representation. Using this result, it was proved in [HOTh] that if $M$ is a compact, full, irreducible isoparametric submanifold of a finite dimensional Euclidean space and $\xi$ is a parallel normal vector field on $M$ such that the focal submanifold $M_\xi$ is not a point, then the fibers of the projection $\pi_\xi : M \longrightarrow M_\xi$, where $\pi_\xi(x) = x + \xi(x)$, are orbits of s-representations. This theorem is called the homogeneous slice theorem and it has been used in both proofs of Thorbergsson's theorem by Thorbergsson and Olmos. In this section, we will generalize the homogeneous slice theorem to finite dimensional slices of infinite dimensional isoparametric submanifolds. Notice that in this case $M_\xi$ is never a point. In fact, it is a proper Fredholm submanifold of the Hilbert space $V$, and therefore always has infinite dimension. Moreover, the connected components of the fibers of the $\pi_\xi$ for various $\xi$ are precisely the finite dimensional $L_P$'s, i.e. those for which $P$ does not contain 0.

We first generalize a part of Olmos' theorem on normal holonomy groups to the infinite dimensional case.

LEMMA 2.1. *Let $N$ be a connected proper Fredholm submanifold of a Hilbert space $V$. For any $x \in N$, let $\Phi(x)$ (respectively $\Phi^*(x)$) be the normal holonomy group (respectively the restricted normal holonomy group) of $N$ at $x$. Then $\Phi(x)$ is a Lie group with $\Phi^*(x)$ as its connected component of the unit element. If all shape operators of $N$ are Hilbert-Schmidt and $\Phi^*(x)$ acts irreducibly on $\nu_x N$, then $\Phi^*(x)$ is a compact Lie group and its representation on $\nu_x N$ is an s-representation or the 1-dimensional trivial representation.*

*Proof.* It follows from Morse theory that the fundamental group of $N$ is countable (cf. [T2], and [PT2, Th. 9.7.6]). Since the codimension of $N$ is finite, the same proof as in the finite dimensional case shows that $\Phi(x)$ is a Lie subgroup of $O(\nu_x N)$ and $\Phi^*(x)$ is its connected component of the unit element (cf. [KN, Ch. 2, §4]). Let $R^\perp$ be the curvature tensor of the normal vector bundle $\nu N$. Then the standard argument shows that $\Phi^*(x) = \{1\}$ if and only if $R^\perp \equiv 0$. Moreover, if $\Phi^*(x)$ acts on $\nu_x N$ irreducibly, then it is a compact subgroup of $O(\nu_x N)$ (cf. [KN, App. 5]).

To prove that the representation of $\Phi^*(x)$ on $\nu_x N$ is an s-representation, we will use the same idea as in the finite dimensional case (cf. [O2]) and emphasize only on those parts where additional arguments are needed. Let $\mathfrak{L}^2$ be the space of Hilbert-Schmidt operators on $T_x N$, and $\langle \cdot, \cdot \rangle_{\mathfrak{L}^2}$ the Hilbert-Schmidt inner product on $\mathfrak{L}^2$, i.e. $\langle B_1, B_2 \rangle_{\mathfrak{L}^2} = \text{Trace}(B_2^* B_1)$. Let $A$ be the shape operator of $N$. For $v_i \in \nu_x N$, $i = 1, 2$, define $\mathfrak{R}^\perp(v_1, v_2) : \nu_x N \longrightarrow \nu_x N$ by

$$\left\langle \mathfrak{R}^\perp(v_1, v_2) v_3, v_4 \right\rangle = (-1/2) \left\langle [A_{v_1}, A_{v_2}], [A_{v_3}, A_{v_4}] \right\rangle_{\mathfrak{L}^2}$$

for all $v_3, v_4 \in \nu_x N$. It is routine to check that $\mathfrak{R}^\perp$ is a curvature tensor on $\nu_x N$. Moreover, by the Ricci equation, for an orthonormal basis $\{e_i \mid i = 1, 2, 3, \cdots\}$ of $T_x N$,

$$\left\langle \mathfrak{R}^\perp(v_1, v_2) v_3, v_4 \right\rangle = \sum_i \left\langle R^\perp(A_{v_1} e_i, A_{v_2} e_i) v_3, v_4 \right\rangle$$

(cf. [O1]). It follows that the triple $(\nu_x N, \mathfrak{R}^\perp, \Phi^*(x))$ is an irreducible holonomy system in the sense of Simons (cf. [S]). If $\mathfrak{R}^\perp$ has nonzero scalar curvature, then by a well known lemma of Simons (cf. [S]), the action of $\Phi^*(x)$ on $\nu_x N$ is an s-representation. On the other hand, if the scalar curvature of $\mathfrak{R}^\perp$ is zero, then $\mathfrak{R}^\perp \equiv 0$ since by the definition, $\mathfrak{R}^\perp$ has nonpositive sectional curvature. This implies that $[A_{v_1}, A_{v_2}] = 0$ for all $v_1, v_2 \in \nu_x N$. By the Ricci equation, $R^\perp \equiv 0$ at $x$. Since the representations of the normal holonomy group on other normal spaces are conjugate to the one on $\nu_x N$, we may also assume that $R^\perp \equiv 0$ at all points of $N$. This is equivalent to saying that $\Phi^*(x) = \{1\}$. The proof is thus finished. □

We now prove the homogeneous slice theorem for finite dimensional slices. We will state the result in terms of horizontal curves which is more convenient for later applications. Let $M$ be a full irreducible isoparametric submanifold of an infinite dimensional Hilbert space $V$. Assume that the rank of $M$ is at least two. Fix $x_0 \in M$ and let $P$ be an affine subspace of $\nu_{x_0} M$. For any horizontal curve $\gamma$ with respect to the distribution $D_P$, let $f_\gamma : W_P(\gamma(0)) \longrightarrow W_P(\gamma(1))$ be the isometry defined in Proposition 1.1. For any $x \in M$, let $\Phi_P(x)$ be the group of isometries of $W_P(x)$ generated by

$$\{f_\gamma \mid \gamma(0), \ \gamma(1) \in L_P(x), \text{ and } \gamma \text{ is horizontal with respect to } D_P\}.$$

If $P$ does not contain 0, then $W_P(x)$ is a finite dimensional Euclidean space and, by the remark following Proposition 1.1, $\Phi_P(x)$ is the normal holonomy group of a focal submanifold. Therefore $\Phi_P(x)$ is a Lie subgroup of the isometry group of $W_P(x)$. Let $\Phi_P^*(x)$ be the connected component of the unit element of $\Phi_P(x)$.

PROPOSITION 2.2. *Let $P$ be an affine subspace of $\nu_{x_0} M$ which does not contain 0. Then $\Phi_P^*(x)$ acts transitively on $L_P(x)$. Assume that $L_P(x)$ is an irreducible isoparametric submanifold of $W_P(x)$. Then the representation of $\Phi_P^*(x)$ on $W_P(x)$ is an irreducible s-representation which has $L_P(x)$ as a principal orbit.*

*Remark.* Although the action of $\Phi_P^*(x)$ on $W_P(x)$ is strictly speaking by affine transformations, we view it as a linear representation. This is justified by the fact that $L_P(x)$ lies in a unique sphere of $W_P(x)$ whose center is of course fixed by $\Phi_P^*(x)$.

*Proof.* We first observe that, since $L_P(x)$ is connected, $\Phi_P^*(x)$ acts on $L_P(x)$ transitively if and only if $\Phi_P(x)$ does. Since $P$ does not contain 0, both $L_P(x)$ and $W_P(x)$ are finite dimensional. By Lemma 3.1 of [HOTh], to prove the transitivity of the $\Phi_P(x)$ action on $L_P(x)$, it suffices to show that for any pair of points $y_1, y_2 \in L_P(x)$ such that $y_2$ lies in some curvature sphere through $y_1$, there exists $f_\gamma \in \Phi_P(x)$ such that $f_\gamma(y_1) = y_2$. This is equivalent to saying that there exists a horizontal curve $\gamma$ with respect to $D_P$ which connects $y_1$ and $y_2$. Assume that there exists a curvature normal $n_i(x_0) \in P$ such that $y_2 \in S_i(y_1)$. Since $M$ is infinite dimensional and irreducible, its Coxeter group is infinite and irreducible (cf. [T2] and [HL]). Therefore there exists a curvature normal $n_j$ such that $n_j(x_0) \notin P$ and $n_j(x_0)$ is neither perpendicular nor parallel to $n_i(x_0)$. Let $l_{i,j}$ be the affine line in $\nu_{x_0} M$ which passes through $n_i(x_0)$ and $n_j(x_0)$. Then $L_{l_{i,j}}(y_1)$, the leaf of distribution $D_{l_{i,j}}$, is a finite dimensional irreducible rank 2 isoparametric submanifold which contains $S_i(y_1)$. By the finite dimensional homogeneous slice theorem, there exists a curve $\gamma$ in $L_{l_{i,j}}(y_1)$ which connects $y_1$ and $y_2$, and which is horizontal with respect to the curvature distribution $E_i$. Since $l_{i,j} \cap P = \{n_i(x_0)\}$, $\gamma$ is also a horizontal curve with respect to $D_P$ when considered as a curve in $M$. This proves that $\Phi_P^*(x)$ acts transitively on $L_P(x)$.

We now prove that the representation of $\Phi_P^*(x)$ on $W_P(x)$ is an $s$-representation. Since $L_P(x)$ is a compact isoparametric submanifold of $W_P(x)$, there is a parallel normal vector field $r_P$ on $M$ such that $x + r_P(x) \in W_P(x)$ and $L_P(x)$ is contained in the round sphere in $W_P(x)$ with center at $x + r_P(x)$ and radius $\|r_P(x)\|$. Let $\pi : M \longrightarrow M_{r_P}$ be the projection of $M$ onto the focal manifold $\bar{M} := M_{r_P}$ which is defined by $\pi(x) = x + r_P(x)$. For any $\bar{x} \in \bar{M}$, let $\bar{W}_P(\bar{x}) = W_P(x) - \bar{x}$, where $x \in \pi^{-1}(\bar{x})$. Also, $\bar{W}_P(\bar{x})$ is a linear subspace of $\nu_{\bar{x}}\bar{M}$ which does not depend on the choice of $x \in \pi^{-1}(\bar{x})$. Let $\bar{\Phi}_P(\bar{x})$ be the normal holonomy group of $\bar{M}$ at $\bar{x}$ which acts on $\nu_{\bar{x}}\bar{M}$. The action of $\bar{\Phi}_P(\bar{x})$ is trivial on the orthogonal complement of $\bar{W}_P(\bar{x})$, and, when restricted to $\bar{W}_P(\bar{x})$, is equivalent to the $\Phi_P(x)$ action on $W_P(x)$. Let $\bar{\Phi}_P^*(\bar{x})$ be the identity component of $\bar{\Phi}_P(\bar{x})$. The action of $\bar{\Phi}_P^*(\bar{x})$ on $\bar{W}_P(\bar{x})$ is irreducible since one of its orbits, $L_P(x)$, is an irreducible isoparametric submanifold of $\bar{W}_P(\bar{x})$. If we can prove that the shape operator of $\bar{M}$ along any vector in $\bar{W}_P(\bar{x})$ is Hilbert-Schmidt, then the same proof as in Lemma 2.1 shows that the representation of $\bar{\Phi}_P^*(\bar{x})$ on $\bar{W}_P(\bar{x})$ is an $s$-representation, and the proposition follows.

For any $x \in M$, the group generated by the reflections along hyperplanes $\{v \in \nu_x M \mid \langle v, n_i(x) \rangle = 1\}$, $i \in I$, is a Coxeter group. Therefore in this set of hyperplanes, there are only a finite number of them which are not parallel to each other, and for each family of parallel ones, the distance between two consecutive hyperplanes is constant. Hence there are only finitely many non-proportional curvature normals, and for each family of proportional ones, the length of the $k^{\text{th}}$ curvature normal is given by the formula $c_1/(c_2+k)$, where $c_1$ and $c_2$ are constants. Moreover, there are only finitely many distinct multiplicities for the curvature distributions. Therefore the series of the eigenvalues of each shape operator of $M$ at $x$ is square summable. Consequently, each shape operator of $M$ is Hilbert-Schmidt. Let $A$ and $\bar{A}$ be the shape operators of $M$ and $\bar{M}$, respectively. For any $x \in M$ and $\bar{x} = \pi(x)$, $T_{\bar{x}}\bar{M}$ is a linear subspace of $T_x M$ which is preserved by the shape operator of $M$ at $x$. Moreover $r_P(x)$ is a normal vector to both $M$ and $\bar{M}$ at $x$ and $\bar{x}$, respectively. The relation between $A$ and $\bar{A}$ is given by the formula:

$$\bar{A}_{r_P(x)} = A_{r_P(x)} \mid_{T_{\bar{x}}\bar{M}} \circ \left( \text{Id} - A_{r_P(x)} \mid_{T_{\bar{x}}\bar{M}} \right)^{-1}.$$

Since $A_{r_P}(x)$ is Hilbert-Schmidt, this formula implies that $\bar{A}_{r_P(x)}$ is also a Hilbert-Schmidt operator on $T_{\bar{x}}\bar{M}$. Since the set $\{r_P(x) \mid x \in \pi^{-1}(\bar{x})\}$ spans $\bar{W}_P(\bar{x})$, $\bar{A}_v$ is Hilbert-Schmidt for all $v \in \bar{W}_P(\bar{x})$. The proof of the proposition is thus finished. $\square$

Let $P_1$ and $P_2$ be two affine subspaces of $\nu_{x_0} M$ such that $P_1 \subset P_2$, $P_1 \neq P_2$ and $0 \notin P_2$. Fix $x \in M$; there are three groups acting on $W_{P_1}(x)$. The first one is $\Phi_{P_1}(x)$, which is abbreviated as $\Phi_1$. By Lemma 1.4, the subgroup of $\Phi_{P_2}(x)$ which is generated by the isometries induced by those curves $\gamma$ which are horizontal with respect to $D_{P_2}$ and $\gamma(0), \gamma(1) \in L_{P_1}(x)$ also acts on $W_{P_1}(x)$. This group is the isotropy group of $\Phi_{P_2}(x)$ at $\bar{x}$, where $\bar{x}$ is the center of the unique sphere of $W_{P_1}(x)$ which contains $L_{P_1}(x)$. Its action on $W_{P_1}(x)$ is not necessarily effective. We denote the corresponding group which acts effectively on $W_{P_1}(x)$ by $(\Phi_2)'_{\bar{x}}$. Since $L_{P_2}(x)$ is an isoparametric submanifold of $W_{P_2}(x)$ and the restriction of $D_{P_1}$ to $L_{P_2}(x)$ is an integrable distribution with totally geodesic leaves, the horizontal curves of $D_{P_1}$ which lie in $L_{P_2}(x)$ also generate a group which acts isometrically on $W_{P_1}(x)$. We denote this group by $\Phi_{1,2}$. More precisely, $\Phi_{1,2}$ is the group of isometries of $W_{P_1}(x)$ generated by

$$\{f_\gamma \mid \gamma(0),\ \gamma(1) \in L_{P_1}(x),\ \gamma(t) \in L_{P_2}(x) \text{ for all } t,$$
$$\text{and } \gamma \text{ is horizontal with respect to } D_{P_1}\}.$$

These three groups have the following important relationship

PROPOSITION 2.3. *If $L_{P_2}(x)$ is irreducible, then $\Phi_1 = (\Phi_2)'_{\bar{x}} = \Phi_{1,2}$ and their representations on $W_{P_1}(x)$ coincide.*

*Remark.* What we will need later is only the equality $\Phi_1 = \Phi_{1,2}$, which can be interpreted also as follows. A connected, complete irreducible isoparametric submanifold $M$ induces on any finite dimensional submanifold of the form $L_P(x)$, which is irreducible when viewed as an isoparametric submanifold of $W_P(x)$ and which is different from $M$ (if $M$ is finite dimensional), the additional structure of a transitive group action. If $L_P(x) \subset L_{P'}(x)$, $L_P(x) \neq L_{P'}(x)$, and if $L_{P'}(x)$ is also finite dimensional and irreducible, then the additional structures on $L_P(x)$ induced from $M$ and $L_{P'}(x)$ coincide. We do not know whether the assumptions $\dim(L_P(x))$, $\dim(L_{P'}(x)) < \infty$ are really necessary.

*Proof.* By the remark following Proposition 1.1, we can restate the above relationship by using normal holonomy groups of focal submanifolds of $M$. The first equality is then just the infinite dimensional version of Theorem 1.2 of [O2] whose proof also applies to infinite dimensional isoparametric submanifolds without any difficulty. By Proposition 2.2, the representation of $\Phi_{P_2}(x)$ on $W_{P_2}(x)$ is an irreducible $s$-representation. Therefore the second equality follows from Theorem 2 of [HO] (cf. also the discussion in 2.2 of [O2]). □

## 3. Determining isoparametric submanifolds by their slices

In this section, we study how much local information is needed to uniquely determine an isoparametric submanifold. Fix $p \in M$. For every affine line $l$ in $\nu_p M$, there is a totally geodesic submanifold $L_l(p)$ of $M$ which is the leaf of the distribution $D_l$ through $p$. If $l$ does not contain any curvature normal, then $L_l(p) = \{p\}$. If $l$ contains exactly one curvature normal, then $L_l(p)$ is either a curvature sphere or $p + E_0(p)$. If $l$ contains more than one curvature normal, then, depending on whether $l$ contains 0 or not, $L_l(p)$ is either a compact rank two isoparametric submanifold or a rank one isoparametric submanifold which contains $p + E_0(p)$. In the latter case, $L_l(p)$ is in general infinite dimensional if $M$ is infinite dimensional. If $M$ is irreducible, then the set of all $L_l(p)$, where $l$ ranges over all affine lines of $\nu_p M$ except one, determines $M$. More precisely, we have:

PROPOSITION 3.1. *Let $M$ and $M'$ be two irreducible isoparametric submanifolds of $V$ with rank bigger than or equal to 2. Assume that there exist $p \in M \cap M'$ such that $T_p M = T_p M'$ and a 1-dimensional linear subspace $l_0$ of $\nu_p M$ satisfying the following condition:*

$(*)$ $\quad L_l(p) = L'_l(p)$ *for any affine line $l \subset \nu_p M$ which is different from $l_0$,*

*where "$'$" denotes the corresponding objects in $M'$. Then $M = M'$.*

*Remark.* For the proof of Thorbergsson's Theorem (i.e. the case $\dim(M) < \infty$), the stronger assumption "$L_l(p) = L'_l(p)$ for all affine lines $l \subset \nu_p M$" would suffice. This shortens somewhat the proof of Lemma 3.2 below. Moreover for a compact finite dimensional isoparametric submanifold $M$, $M = Q(p)$ for any $p \in M$. Therefore the arguments in Case (1) of the proof of Lemma 3.2 would be sufficient.

To prove Proposition 3.1 for the infinite dimensional case, we need to introduce another equivalence relation. For any two points $x$ and $y$ in $M$, we say $x \sim_0 y$ if there exists a sequence of points $x_k \in M$, $k = 0, 1, \cdots, n$, such that $x_0 = x$, $x_n = y$, and for every $k = 1, \cdots, n$, either $x_k \sim x_{k-1}$ or $x_k \in x_{k-1} + E_0(x_{k-1})$. Recall that the equivalence relation "$\sim$" was defined in the introduction. For any $x \in M$, we define

$$Q_0(x) := \{y \in M \mid y \sim_0 x\}.$$

We first prove the following

LEMMA 3.2. *If $M$ and $M'$ satisfy the conditions in Proposition 3.1, then $\overline{Q_0(p)} = \overline{Q'_0(p)}$, where $\overline{Q_0(p)}$ is the closure of $Q_0(p)$ and $Q'_0(p)$ is the set of all $q \in M'$ such that $q \sim_0 p$ in $M'$. Moreover condition $(*)$ holds for all points in $\overline{Q_0(p)}$.*

*Proof.* Since both $M$ and $M'$ are complete, we only need to show that $Q_0(p) = Q'_0(p)$ and condition $(*)$ holds for all points in $Q_0(p)$. For every $q \in Q_0(p)$, there exists a finite sequence of points $x_k$, where $k = 0, 1, \ldots, m$, such that $x_0 = p$, $x_m = q$, and $x_k$ lies either in $x_{k-1} + E_0(x_{k-1})$ or in some curvature sphere of $x_{k-1}$ for $k = 1, \ldots, m$. Suppose that for some $k$, $x_k$ lies in a curvature sphere $S_r(x_{k-1})$ whose curvature normal $n_r \in l_0$. Since $M$ is irreducible and has rank at least two, there exists a curvature normal $n_s$ which is neither parallel nor perpendicular to $n_r$. Let $l_{r,s}$ be the affine line in $\nu_{x_{k-1}}M$ which passes through $n_r(x_{k-1})$ and $n_s(x_{k-1})$. Then the leaf $L_{l_{r,s}}(x_{k-1})$ is a finite dimensional irreducible rank 2 isoparametric submanifold which contains $S_r(x_{k-1})$. By Theorem D of [HOTh], $x_{k-1}$ can be joined to $x_k$ by a piecewise differentiable curve whose pieces are tangent to one of $E_j$ where $n_j \in l_{r,s}$ and $n_j \neq n_r$. Therefore, enlarging the sequence $\{x_k \mid k = 0, 1, \ldots, m\}$ if necessary, we can assume that for $k = 1, \ldots, m$, $x_k$ lies either in $x_{k-1} + E_0(x_{k-1})$ or in a curvature sphere of $x_{k-1}$ whose corresponding curvature normal does not belong to $l_0$. Hence, to prove that $q \in M'$ and satisfies the condition $(*)$, it suffices to prove that condition $(*)$ holds for every point $q$ which lies either in $p + E_0(p)$ or in some curvature sphere of $p$ whose corresponding curvature normal does not belong to $l_0$. The lemma then follows from a trivial induction argument on $m$.

*Case* 1. $q \in S_k(p)$, where $n_k \notin l_0$.

By condition $(*)$ at $p$, $S_k(q) = S_k(p) = S'_k(p) = S'_k(q)$. For any affine line $l \subset \nu_p M$ such that $n_k \in l$, $q \in S_k(p) \subset L_l(p) \cap L'_l(p)$. Since $l \neq l_0$, by condition $(*)$ at $p$, $L_l(q) = L_l(p) = L'_l(p) = L'_l(q)$. This implies that for any $n_j \in l$, $S_j(q) = S'_j(q)$ if $n_j \neq 0$, and $E_0(q) = E'_0(q)$ if $0 \in l$. Therefore $S_j(q) = S'_j(q)$ for all $j$ because for every curvature normal $n_j \neq n_k$, there is a unique affine line passing through $n_j$ and $n_k$. In particular, this implies that at $q$, the curvature normals of $M$ and $M'$ coincide.

To show that condition $(*)$ holds at $q$, it remains to show that for any affine line $l \subset \nu_p M$ such that $n_k \notin l$ and $l \neq l_0$, $L_l(q) = L'_l(q)$. Choose a smooth curve $\gamma$ in $S_k(p)$ such that $\gamma(0) = p$ and $\gamma(1) = q$. Then $\gamma$ is a horizontal curve with respect to the distribution $D_l$. Therefore it induces a one-parameter family of isometries $f_t : W_l(p) \longrightarrow W_l(\gamma(t))$ by Proposition 1.1. On the other hand, as a horizontal curve in $M'$, $\gamma$ also induces a one-parameter family of isometries $f'_t : W'_l(p) \longrightarrow W'_l(\gamma(t))$, where $W'_l$ is defined in the same way as for $W_l$. By Proposition 1.1, $f_t(L_l(p)) = L_l(\gamma(t))$ and $f'_t(L'_l(p)) = L'_l(\gamma(t))$. Since $p$ satisfies condition $(*)$, $L_l(p) = L'_l(p)$ and $W_l(p) = W'_l(p)$. Therefore to prove that $L_l(\gamma(t)) = L'_l(\gamma(t))$, it suffices to show that $f_t \equiv f'_t$ for all $t$. Since both $f_t$ and $f'_t$ are isometries, we only need to show that $(f_t)_* = (f'_t)_*$. Since the curvature normals of $M$ and $M'$ at both $p$ and $q$ coincide with each other, the parallel translation in $\nu M$ from $p$ to $q$ coincides with that in $\nu M'$. Therefore, by Lemma 1.3, $(f_t)_*$ and $(f'_t)_*$ are equal when restricted to the normal space of $L_l(p)$ in $W_l(p)$. For any curvature normal $n_r(p) \in l$, let $l_{r,k}$ be the affine line in $\nu_p M$ which passes through $n_r(p)$ and $n_k(p)$. Since $\gamma$ is a horizontal curve in $L_{l_{r,k}}(p) = L'_{l_{r,k}}(p)$ with respect to $E_r \mid_{L_{l_{r,k}}(p)} = E'_r \mid_{L'_{l_{r,k}}(p)}$, it induces two one-parameter families of isometries, denoted by $g_t$ and $g'_t$, in $L_{l_{r,k}}(p)$ and $L'_{l_{r,k}}(p)$ respectively. We have $g_t = g'_t$. By Lemma 1.4, $(f_t)_* |_{E_r(p)} = (g_t)_* |_{E_r(p)} = (g'_t)_* |_{E'_r(p)} = (f'_t)_* |_{E'_r(p)}$. This proves that $(f_t)_* = (f'_t)_*$.

*Case 2.* $q \in p + E_0(p)$.

For any affine line $l \subset \nu_p M$ such that $0 \in l$ and $l \neq l_0$, we have $L_l(q) = L_l(p) = L'_l(p) = L'_l(q)$. In particular, for any curvature normal $n_k \notin l_0$, we can take $l$ to be the affine line which passes through $0$ and $n_k(p)$. Let $\gamma$ be the straight line in $p + E_0(p)$ which connects $p$ to $q$. By Proposition 1.1, $\gamma$ induces two one-parameter families of isometries $f_t : S_k(p) \longrightarrow S_k(\gamma(t))$ and $f'_t : S'_k(p) \longrightarrow S'_k(\gamma(t))$ when it is considered as a horizontal curve with respect to $E_k$ and $E'_k$ respectively. The above relation implies that $f_t = f'_t$ for all $t$.

It remains to show that if $l$ is an affine line in $\nu_p M$ such that $0 \notin l$ and $l \neq l_0$, then $L_l(q) = L'_l(q)$. In this case the straight line $\gamma$ defined above is a horizontal curve with respect to the distribution $D_l$ on $M$ as well as to the distribution $D'_l$ on $M'$. Therefore it induces two one-parameter families of isometries $f_t : L_l(p) \longrightarrow L_l(\gamma(t))$ and $f'_t : L'_l(p) \longrightarrow L'_l(\gamma(t))$. Since $L_l(p) = L'_l(p)$, it suffices to show that $(f_t)_* |_{D_l(p)} = (f'_t)_* |_{D'_l(p)}$ for all $t$. By Lemma 1.4 and the arguments in the previous paragraph, $(f_t)_* |_{E_k(p)} = (f'_t)_* |_{E'_k(p)}$ for all $n_k(p) \in l$ such that $n_k(p) \notin l_0$. Since $D_l(p) = \bigoplus_{n_k(p) \in l} E_k(p)$, we are done if $l \cap l_0$ does not contain any curvature normal. Therefore we can assume that $l \cap l_0 = \{n_i(p)\}$ for some nonzero curvature normal $n_i(p)$. We only need to show that $(f_t)_* |_{E_i(p)} = (f'_t)_* |_{E'_i(p)}$.

If $L_l(p)$ is an irreducible isoparametric submanifold of rank 2, then there exists a curvature normal $n_j(p) \in l$ which is neither perpendicular nor parallel to $n_i(p)$. Let $\tilde{p}$ be the antipodal point of $p$ in the curvature sphere $S_j(p)$. Then $\tilde{p} \in L_l(p)$ and $E_i(p) = E_k(\tilde{p})$ for some curvature distribution $E_k$ with curvature normal $n_k(p) \notin l_0$ (cf. [T2, Th. 6.11 (i) and (iii)]). Let $\tilde{\gamma}(t) = \gamma(t) + 2n_j(\gamma(t))/\|n_j\|^2$. Notice that $n_j(\gamma(t))$ is constant since $\gamma(t) \in p + E_0(p)$ for all $t$. Hence $\tilde{\gamma}$ is a straight line in $\tilde{p} + E_0(\tilde{p})$ (cf. [T2, Th. 6.11 (v)]). It follows that $\tilde{\gamma}$ is a horizontal curve with respect to $D_l$ which is parallel to $\gamma$. Moreover $\tilde{\gamma}(t) \in M'$ for all $t$ since $M$ and $M'$ have same curvature normals at $p$ and $E_0(p) = E'_0(p)$. Therefore $f_t$ and $f'_t$ can also be viewed as one-parameter families of isometries induced by $\tilde{\gamma}$. By the result from Case (1), we know that condition $(*)$ also holds for $\tilde{p}$. Therefore, as in the previous paragraph, we have that at the point $\tilde{p}$, $(f_t)_* |_{E_k(\tilde{p})} = (f'_t)_* |_{E'_k(\tilde{p})}$. Since both $f_t$ and $f'_t$ are restrictions of affine maps from $W_l(p)$ to $V$, this implies that at the point $p$, $(f_t)_* |_{E_i(p)} = (f'_t)_* |_{E'_i(p)}$.

On the other hand, since the Coxeter group of $M$ is irreducible (cf. [HL]), for each nonzero curvature normal $n_i(p) \in l_0$, there exists a curvature normal $n_j(p)$ which is neither perpendicular nor parallel to $n_i(p)$. Let $l_{i,j}$ be the affine line in $\nu_p(M)$ passing through $n_i(p)$ and $n_j(p)$. Then the corresponding slice $L_{l_{i,j}}(p)$ is a finite dimensional irreducible homogeneous rank 2 isoparametric submanifold. Using the above argument for this slice, we have that $E_i(\gamma(t)) = E'_i(\gamma(t))$ for all $t$. This proves that $L_l(\gamma(t)) = L'_l(\gamma(t))$ for all $t$ if $l$ contains only one curvature normal $n_i(p)$ or if $L_l(p)$ is a reducible rank 2 isoparametric submanifold. The proof of the lemma is thus finished. □

*Proof of Proposition* 3.1. By Theorem D, $Q_0(x)$ is dense in $M$. Therefore the proposition follows from Lemma 3.2. □

## 4. Constructing isometries

Assume that $M$ is a complete, irreducible, full isoparametric submanifold in a Hilbert space $V$ and the set of all the curvature normals of $M$ at one point is not contained in any affine line. The dimension of $M$ could be either finite or infinite. As mentioned in the introduction, the last condition is equivalent to saying that $M$ has codimension at least 3 in the finite dimensional case and at least 2 in the infinite dimensional case. For any $p \in M$ and curvature normal $n_i \neq 0$, let $\Phi_i(p)$ be the group of isometries of $W_i(p) := p + \mathbb{R}n_i(p) \oplus E_i(p)$ generated by

$$\{f_\gamma \mid \gamma(0),\ \gamma(1) \in S_i(p),\ \text{and}\ \gamma\ \text{is horizontal with respect to}\ E_i\}.$$

Let $\Phi_i^*(p)$ be the identity component of $\Phi_i(p)$. To construct isometries on $V$ which preserve $M$, we need an appropriate biinvariant metric on $\Phi_i^*(p)$. Let $\rho : K \longrightarrow O(\mathbb{E})$ be an orthogonal representation of a compact Lie group $K$ on a finite dimensional Euclidean vector space $\mathbb{E}$. Then we define a biinvariant metric on $K$ by

$$\langle X, Y \rangle := -B(X, Y) - \text{Trace}(\rho_* X)(\rho_* Y)$$

for any $X$ and $Y$ from the Lie algebra of $K$, where $B$ denotes the Killing form. We call this the metric on $K$ associated to $\rho$. In case of the isotropy representation of a symmetric space $G/K$, the corresponding metric on $K$ is the metric induced from the Killing form of the Lie algebra of $G$. We equip $\Phi_i^*(p)$ with the metric associated to its representation on $W_i(p)$, which can be viewed as a linear representation since it fixes the center of $S_i(p)$. For any $X$ in the Lie algebra of $\Phi_i^*(p)$ which is perpendicular to the Lie algebra of the isotropy subgroup of $\Phi_i^*(p)$ at $p$, define $g(t) = \exp(tX) \in \Phi_i^*(p)$ and $\gamma(t) = g(t) \cdot p$. Then $\gamma$ is a curve in $S_i(p)$ with $\gamma(0) = p$. Let $q = \gamma(1)$. For any curvature normal $n_j$ which is different from $n_i$ (here $n_j$ could be zero), $\gamma$ is a horizontal curve with respect to $E_j$. Let $f_\gamma^j$ be the isometry from $W_j(p) = p + \mathbb{R}n_j(p) \oplus E_j(p)$ to $W_j(q)$ defined by Proposition 1.1. Then $f_\gamma^j$ maps $S_j(p)$ to $S_j(q)$. Let $F_\gamma$ be the unique isometry of the ambient Hilbert space $V$ which satisfies the following properties: $F_\gamma(p) = q$, $(F_\gamma)_* |_{\nu_p M}$ is the parallel translation in $\nu M$ from $p$ to $q$, $(F_\gamma)_* |_{E_i(p)} = g(1)_* |_{E_i(p)}$, and for any $n_j$ which is different from $n_i$, $(F_\gamma)_* |_{E_j(p)} = (f_\gamma^j)_* |_{E_j(p)}$. The main result of this section is the following:

THEOREM 4.1. *$F_\gamma(M) = M$, $(F_\gamma)_*$ preserves each curvature distribution on $M$, and $(F_\gamma)_* |_{\nu M}$ coincides with the parallel translation in the normal bundle.*

We begin with a special case where $M$ is a finite dimensional homogeneous isoparametric submanifold. Notice that the definition of $F_\gamma$ does not use the restriction on the codimension of $M$. In the following lemma, the codimension of $M$ could also be 2. In fact, we will only need that case later.

LEMMA 4.2. *If $M$ is the principal orbit of an irreducible s-representation of any codimension (and hence an irreducible compact homogeneous isoparametric submanifold), then $F_\gamma(M) = M$.*

*Proof.* Let $G/K$ be an irreducible symmetric space of noncompact type where $G$ is the identity component of the group of isometries of $G/K$ and $\mathfrak{g} = \mathfrak{k} + \mathfrak{p}$ the corresponding Cartan decomposition of the Lie algebra of $G$. Let $B$ be the Killing form of $\mathfrak{g}$. Then $B$ is positive definite on $\mathfrak{p}$ and negative definite on $\mathfrak{k}$. $B$ and $-B$ induce a well defined $\mathrm{Ad}(K)$-invariant inner product on $\mathfrak{p}$ and $\mathfrak{k}$ respectively. Assume that $M = \mathrm{Ad}(K)v \subset \mathfrak{p}$ for some $v \in \mathfrak{p}$ which is a regular point of the $\mathrm{Ad}(K)$-action. Let $\mathfrak{a}$ be a maximal abelian subspace of $\mathfrak{p}$ which contains $v$ and let $\Lambda \subset \mathfrak{a}^* \setminus \{0\}$ be the set of roots with respect to $\mathfrak{a}$ (cf. [Lo]). Since $M$ is a principal orbit of the $\mathrm{Ad}(K)$-action, $\mathfrak{a} = \nu_v M$ and $\lambda(v) \neq 0$ for all $\lambda \in \Lambda$. For any fixed curvature normal $n_i$ of $M$, let $r_i = n_i/\|n_i\|^2$. Considering parallel submanifolds if necessary, we may assume that for any curvature normal $n_j$, $\langle n_j, r_i \rangle = 1$ if and only if $n_j = n_i$. Let $v_i = v + r_i(v) \in \mathfrak{a}$. Then $v_i$ is a singular point of the $\mathrm{Ad}(K)$-action and $M_i := \mathrm{Ad}(K)v_i$ is a focal submanifold of $M$ with the curvature sphere $S_i(v)$ as a slice. Let $\Lambda_i^0 = \{\lambda \in \Lambda \mid \lambda(v_i) = 0\}$ and $\Lambda_i^+ = \{\lambda \in \Lambda \mid \lambda(v_i) > 0\}$. Let $\mathfrak{k}_0$ be the centralizer of $\mathfrak{a}$ in $\mathfrak{k}$ and for any $\lambda \in \Lambda$, let $\mathfrak{k}_\lambda = \{X \in \mathfrak{k} \mid (\mathrm{ad}A)^2 X = \lambda^2(A)X \text{ for all } A \in \mathfrak{a}\}$ and $\mathfrak{p}_\lambda = \{X \in \mathfrak{p} \mid (\mathrm{ad}A)^2 X = \lambda^2(A)X \text{ for all } A \in \mathfrak{a}\}$. Then there is a decomposition

$$\begin{aligned}\mathfrak{g} &= \mathfrak{k}_0 + \sum_{\lambda \in \Lambda_i^0} \mathfrak{k}_\lambda + \sum_{\lambda \in \Lambda_i^+} \mathfrak{k}_\lambda + \mathfrak{a} + \sum_{\lambda \in \Lambda_i^0} \mathfrak{p}_\lambda + \sum_{\lambda \in \Lambda_i^+} \mathfrak{p}_\lambda \\ &= \mathfrak{k}_{v_i} + \mathfrak{k}_i^+ + \mathfrak{p}_{v_i} + \mathfrak{p}_i^+,\end{aligned}$$

where $\mathfrak{k}_{v_i} = \mathfrak{k}_0 + \sum_{\lambda \in \Lambda_i^0} \mathfrak{k}_\lambda$, $\mathfrak{k}_i^+ = \sum_{\lambda \in \Lambda_i^+} \mathfrak{k}_\lambda$, $\mathfrak{p}_{v_i} = \mathfrak{a} + \sum_{\lambda \in \Lambda_i^0} \mathfrak{p}_\lambda$, $\mathfrak{p}_i^+ = \sum_{\lambda \in \Lambda_i^+} \mathfrak{p}_\lambda$. The geometric interpretation of this splitting is the following (cf. [HO]). $\mathfrak{p}_i^+ = T_{v_i} M_i$ and $\mathfrak{p}_{v_i} = \nu_{v_i} M_i$. The $i^{\mathrm{th}}$ curvature distribution of $M$ at $v$ is $E_i(v) = \sum_{\lambda \in \Lambda_i^0} \mathfrak{p}_\lambda$; $\mathfrak{k}_{v_i}$ is the Lie algebra of $K_{v_i}$, the isotropy group of $K$ at $v_i$; and $\mathfrak{k}_0$ is the Lie algebra of $K_v$, the isotropy group of $K$ at $v$. Moreover $\mathrm{Ad}(K_{v_i})v = S_i(v)$ and $K_v \subset K_{v_i}$ is also the isotropy group of $K_{v_i}$ at $v$. By [HO, Th. 2], the (effectively made) action of $K_{v_i}$ on $\nu_{v_i} M_i$ is equivalent to the action of the normal holonomy group of $M_i$ at $v_i$, and therefore is also equivalent to the action of $\Phi_i(v)$ by the remark following Proposition 1.1. More precisely, let $\mathfrak{k}_0'$ be the orthogonal complement to $\{X \in \mathfrak{k}_0 \mid \mathrm{ad} X = 0 \text{ on } \mathfrak{p}_{v_i}\}$ in $\mathfrak{k}_0$ with respect to the Killing form of $\mathfrak{g}$, and let $\mathfrak{a}'$ be the linear subspace of $\mathfrak{a}$ which is spanned by the root vectors corresponding to elements of $\Lambda_i^0$. Then $\mathfrak{k}_0' + \sum_{\lambda \in \Lambda_i^0} \mathfrak{k}_\lambda$ is the Lie algebra of the normal holonomy group of $M_i$ at $v_i$, and its action on the normal space $\nu_{v_i} M$ leaves invariant the subspace $\mathfrak{a}' + \sum_{\lambda \in \Lambda_i^0} \mathfrak{p}_\lambda$ (this is the subspace corresponding to $W_i(v)$). A straightforward computation shows that $\sum_{\lambda \in \Lambda_i^0} \mathfrak{k}_\lambda$ is the orthogonal complement in the Lie algebra of the normal holonomy group of $M_i$ at $v_i$ to the Lie algebra of the isotropy group at $v$ with respect to the biinvariant metrics associated to the representation on the normal space as well as to the representation on the subspace $\mathfrak{a}' + \sum_{\lambda \in \Lambda_i^0} \mathfrak{p}_\lambda$. The second metric is what we used in the definition of $F_\gamma$.

For any $X \in \sum_{\lambda \in \Lambda_i^0} \mathfrak{k}_\lambda$, let $\gamma(t) = \mathrm{Ad}(\exp(tX))v$ and consider the isometry $F_\gamma$ of $\mathfrak{p}$ defined above. To prove that $F_\gamma(M) = M$, it suffices to show that $F_\gamma = \mathrm{Ad}(\exp(X)) \in \mathrm{Ad}(K)$. Since both $F_\gamma$ and $\mathrm{Ad}(\exp(X))$ are isometries on $\mathfrak{p}$, we only need to show that $(F_\gamma)_* \mid_v = \mathrm{Ad}(\exp(X))_* \mid_v$. By the definition of $F_\gamma$, it follows immediately that $(F_\gamma)_* \mid_{E_i(v)} = \mathrm{Ad}(\exp(X))_* \mid_{E_i(v)}$. Since the $\mathrm{Ad}(K)$-action is polar and $M$ is a principal orbit, $\mathrm{Ad}(\exp(X))_* \mid_{\nu_v M}$ is the parallel translation in $\nu M$, and therefore is equal to $(F_\gamma)_* \mid_{\nu_v M}$. It remains to show that, for any curvature normal $n_j$ which is not equal to $n_i$, $\mathrm{Ad}(\exp(X))_* \mid_{E_j(v)}$ comes from the isometry induced by $\gamma$ as a horizontal curve with respect to $E_j$. It suffices to show that for any $u \in S_j(v)$, the map $t \longmapsto \mathrm{Ad}(\exp(tX))u$ defines a horizontal curve with respect to $E_j$ which is parallel to $\gamma$. Let $v_j = v + n_j(v)/\|n_j(v)\|^2$ and $M_j = \mathrm{Ad}(K)v_j$. By the remark following Proposition 1.1, we only need to show that $\mathrm{Ad}(\exp(tX))$ applied to the normal vector of $M_j$ at $v_j$ defines parallel translation in $\nu M_j$. This follows from [HO, para. 4, p. 872], since $X$ is perpendicular to the Lie algebra of the isotropy group of $K$ at $v_j$ with respect to the Killing form of $\mathfrak{g}$. The proof of the lemma is thus finished. $\square$

Now we return to the general case.

*Proof of Theorem* 4.1. Let $M' = F_\gamma(M)$. Since $F_\gamma$ is an isometry of $V$, $M'$ is also an irreducible isoparametric submanifold of $V$ with the same codimension. By the definition of $F_\gamma$, $q = F_\gamma(p) \in M \cap M'$, $T_q M = T_q M'$, and the curvature normals of $M$ and $M'$ at $q$ coincide with each other. It suffices to show that $q$ satisfies condition $(*)$ of Proposition 3.1 with $l_0$ equal to the line passing through 0 and $n_i(q)$. In fact, Proposition 3.1 will imply $F_\gamma(M) = M$, and the rest of Theorem 4.1 will follow from the last assertion of Lemma 3.2.

For any curvature distribution $E_j$, $E_j(q) = (F_\gamma)_*(E_j(p)) = E'_j(q)$. Therefore all the curvature spheres of $M$ and $M'$ at $q$ coincide. Suppose that $l$ is an arbitrary affine line in $\nu_q M$. If $n_i \in l$, $0 \notin l$, and $l$ contains more than one curvature normal, then $L_l(q)$ is an orbit of the $s$-representation of a rank 2 symmetric space by the homogeneous slice theorem in [HOTh] and Proposition 2.2. Notice that here we used the condition that the set of all curvature normals of $M$ is not contained in any affine line. Moreover $S_i(p) = S_i(q) \subset L_l(q) = L_l(p)$ since $n_i(q) \in l$. If $L_l(q)$ is reducible, then it is a product of $S_i(q)$ with some curvature sphere, say $S_j(q)$. Since $L'_l(q)$ has the same Coxeter group as $L_l(q)$, $L'(q) = S'_i(q) \times S'_j(q)$ (cf. [T1]). Therefore $L_l(q) = L'_l(q)$ because $M$ and $M'$ have the same set of curvature spheres at $q$. If $L_l(q)$ is irreducible, then by Lemma 4.2 (as well as Proposition 2.3), $L_l(q) = F_\gamma(L_l(q)) = F_\gamma(L_l(p)) = L'_l(q)$. If $n_i \notin l$, then $\gamma$ is a horizontal curve with respect to $D_l$. Therefore by Proposition 1.1, $\gamma$ induces a one-parameter family of isometries, denoted by $f_t^l$, from $W_l(p)$ to $W_l(\gamma(t))$. Lemma 1.4 implies that $(f_1^l)_* = (F_\gamma)_* |_{T_p W_l(p)}$. Therefore $f_1^l = F_\gamma |_{W_l(p)}$ since both of them are isometries. Hence $L'_l(q) = F_\gamma(L_l(p)) = f_1^l(L_l(p)) = L_l(q)$ by Proposition 1.1. This proves that $q$ satisfies condition $(*)$. The proof of the theorem is thus finished. □

An immediate consequence of Theorem 4.1 is the homogeneity of $Q(p)$ for $p \in M$ (which already implies Thorbergsson's theorem in the finite dimensional case).

COROLLARY 4.3. *For every $q \in Q(p)$, there exists an isometry $f$ of $V$ such that $f(p) = q$, $f(M) = M$, and $f$ preserves $Q(p)$ and $\overline{Q(p)}$, and for every $x \in M$, $f_* |_{\nu_x M}$ is the parallel translation in $\nu M$ and $f_*(E_i(x)) = E_i(f(x))$ for all $i \in I$.*

*Proof.* It suffices to prove this lemma for the case where $q$ lies in some curvature sphere of $p$. The general case is proved by using the composition of isometries constructed for this special case. Suppose $q \in S_i(p)$. By Proposition 2.2, the group $\Phi_i^*(p)$ acts transitively on $S_i(p)$. Therefore there exists a sequence of points $\{x_k \mid k = 0, 1, \ldots, m\}$ in $S_i(p)$ which satisfies the conditions that $x_0 = p$, $x_m = q$, and for every $k = 1, \ldots, m$, there exists $X_k$ in the Lie algebra of $\Phi_i^*(p)$ such that $X_k$ is perpendicular to the Lie algebra of the isotropy subgroup of $\Phi_i^*(p)$ at $x_{k-1}$ and $(\exp X_k) \cdot x_{k-1} = x_k$. By Theorem 4.1, there exists an isometry, denoted by $F_k$, of $V$ such that $F_k(x_{k-1}) = x_k$, $F_k(M) = M$, and $F_k$ preserves curvature normals and curvature distributions at all points of $M$. Since $F_k$ maps curvature spheres to curvature spheres, it preserves $Q(p)$ and therefore also $\overline{Q(p)}$. Let $f = F_m \circ \cdots \circ F_1$. Then $f$ satisfies the desired properties. □

Now we are ready to prove the homogeneity of the closure of each equivalence class.

PROPOSITION 4.4. *$\overline{Q(p)}$ is extrinsically homogeneous in $V$.*

*Proof.* For any $q \in \overline{Q(p)}$, we want to show that there exists an isometry of $V$ which maps $p$ to $q$ and preserves $\overline{Q(p)}$. Choose a sequence of points $q_k \in Q(p)$, $k = 1, 2, 3, \ldots$, which converges to $q$. By Corollary 4.2, there exists a sequence of isometries, say $f_k$, $k = 1, 2, 3, \ldots$, of $V$ such that $f_k(p) = q_k$, $f_k(\overline{Q(p)}) = \overline{Q(p)}$, and $f_k$ preserves curvature normals and curvature distributions at all points of $\overline{Q(p)}$. It suffices to show that there is a subsequence $\{f_{m_k}\}$ and an isometry $f$ of $V$ such that $\{f_{m_k}(x)\}$ converges to $f$ for all $x$ in the affine span of $\overline{Q(p)}$.

We first notice that if for some $x \in \overline{Q(p)}$, $\{f_k(x)\}$ converges, then for each $i \in I \setminus \{0\}$, there exists a subsequence of $\{f_k\}$ which is pointwise convergent on $S_i(x)$. In fact, since $f_k$ preserves curvature distributions and curvature normals, for every $i \in I \setminus \{0\}$, $f_k(S_i(x)) = S_i(f_k(x))$, which converges to the $i^{\text{th}}$ curvature sphere at the limit of $\{f_k(x)\}$ by the continuity of curvature distributions. Since each $f_k$ is an isometry and $S_i(x)$ is finite dimensional, there exists a subsequence of $\{f_k\}$ which is pointwise convergent on $S_i(x)$.

Since $\{f_k(p)\}$ converges and every point $x \in Q(p)$ can be connected to $p$ by finitely many curvature spheres, using the above fact repeatedly, we can find a subsequence of $\{f_k\}$ which converges at $x$. Notice that this subsequence depends on $x$ in general. Since at each point, there are only countably many curvature spheres, there exists a countable set $C \subset Q(p)$ which is dense in $Q(p)$. Using a standard diagonal process, we can find a subsequence $\{f_{m_k}\}$ of $\{f_k\}$ which is pointwise convergent on $C$. Since each $f_{m_k}$ is an isometry, $\{f_{m_k}\}$ is pointwise convergent on the closure of $C$, which is equal to $\overline{Q(p)}$. Since each $f_{m_k}$ is an affine map, $\{f_{m_k}\}$ is also pointwise convergent on the affine span of $\overline{Q(p)}$. Let $f$ be an arbitrary isometry of $V$ whose restriction on the affine span of $\overline{Q(p)}$ is the pointwise limit of $\{f_{m_k}\}$. Then $f(\overline{Q(p)}) = \overline{Q(p)}$ and $f(p) = q$. The proof is completed. □

We end this section with one more property of the isometry $F_\gamma$ constructed in the beginning of this section. In fact, for each $X$ in the Lie algebra of $\Phi_i^*(p)$ which is perpendicular to the Lie algebra of the isotropy subgroup of $\Phi_i^*$ at $p$, we can define a one-parameter family of isometries $F_X^t$ on $V$ in an obvious way such that $F_X^t(p) = \gamma(t) := \exp(tX) \cdot p$ and $F_X^1 = F_\gamma$.

LEMMA 4.5. *$F_X^t$, $t \in \mathbb{R}$, is a one-parameter group of isometries on $V$ which preserve $M$.*

*Proof.* We only need to show that for all positive numbers $s_1$ and $s_2$, $F_X^{s_2} \circ F_X^{s_1} = F_X^{s_1+s_2}$. Since $F_X^t \mid_{S_i(p)} = \exp(tX)$, we have

$$F_X^{s_2} \circ F_X^{s_1}(p) = F_X^{s_1+s_2}(p) = \gamma(s_1 + s_2),$$

and

$$(F_X^{s_2})_* \circ (F_X^{s_1})_* \mid_{E_i(p)} = (F_X^{s_1+s_2})_* \mid_{E_i(p)} = (\exp(s_1 + s_2))_* \mid_{E_i(p)}.$$

Moreover $(F_X^{s_2})_* \circ (F_X^{s_1})_* \mid_{\nu_p M} = (F_X^{s_1+s_2})_* \mid_{\nu_p M}$ since they are equal to the parallel translation in the normal bundle of $M$. Therefore we only need to show that $(F_X^{s_2})_* \circ (F_X^{s_1})_* \mid_{E_j(p)} = (F_X^{s_1+s_2})_* \mid_{E_j(p)}$ for all $j \neq i$. Let $\tilde{\gamma}(t) = \gamma(s_1 + t)$ and $\tilde{f}_t : S_j(\tilde{\gamma}(0)) \longrightarrow S_j(\tilde{\gamma}(t))$ be the one-parameter family of isometries defined by $\tilde{\gamma}$ as a horizontal curve with respect to $E_j$. It suffices to show that $F_X^{s_2} \mid_{S_j(\tilde{\gamma}(0))} = \tilde{f}_{s_2}$.

Let $c$ be an arbitrary horizontal curve with respect to $E_j$ which starts from a point in $S_j(p)$ and is parallel to $\gamma$. For any $s \geq 0$, $F_X^s \mid_{S_j(p)}$ coincides with the one-parameter family of isometries induced by $\gamma$ as a horizontal curve with respect to $E_j$, because they have the same domains and images and their tangential maps coincide at $p$. Therefore $F_X^s(c(0)) = c(s)$. Since $F_X^s$ is an isometry which preserves the curvature distribution $E_j$, $F_X^s \circ c$ is also a horizontal curve with respect to $E_j$. Moreover it is parallel to $F_X^s \circ \gamma$. On the other hand, since $c$ is parallel to $\gamma$ and $F_X^s \circ \gamma(t) = \gamma(t + s)$, the curve $t \longmapsto c(t + s)$ is also parallel to $F_X^s \circ \gamma$. By the uniqueness of horizontal curves passing through the point $c(s)$, we have $F_X^s(c(t)) = c(s + t)$. In particular, $F_X^{s_2}(c(s_1)) = c(s_2 + s_1) = \tilde{f}(s_2)(c(s_1))$. Since for any $x \in S_j(\tilde{\gamma}(0))$ there always exists a horizontal curve $c$ such that $c(0) \in S_j(p)$, $c(s_1) = x$ and $c$ is parallel to $\gamma$, the proof of the lemma is thus finished. □

## 5. Proof of Theorems A, B, and C

In this section, we prove the theorems in the introduction. Since the finite dimensional case of Theorem A was already settled in Corollary 4.2, we will only deal with infinite dimensional isoparametric submanifolds in this section. We will prove Theorem B first. Let $M$ be an irreducible isoparametric submanifold in a Hilbert space $V$ with codimension at least 2. We call a vector field $w$ on $M$ an *iterated bracket* of vector fields $X_i$, $i = 1, \cdots, k$, if $w = [X_1, [X_2, \cdots [X_{k-1}, X_k] \cdots ]]$. An iterated bracket of this form is of *length* $k - 1$. In the special case where $k = 1$, this iterated bracket should be understood as just $X_1$. For any $x \in M$, let $E'(x)$ be the linear subspace of $T_x M$ consisting of all finite linear combinations of iterated brackets of the form $[X_1, [X_2, \cdots [X_{k-1}, X_k] \cdots ]](x)$, where $k \geq 1$ and each $X_j$, $j = 1, \cdots, k$, lies in a curvature distribution with nonzero curvature normal. By definition, each curvature distribution with nonzero curvature normal is a subspace of $E'(x)$. Let $\overline{E'(x)}$ be the closure of $E'(x)$ and $E''(x)$ be the orthogonal complement of $\overline{E'(x)}$ in $T_xM$. Then $E''(x)$ is a closed subspace of $E_0(x)$. Although it is easy to see that $\overline{E'}$ and $E''$ are distributions on $M$ when $M$ is extrinsically homogeneous, we do not know whether this is true for general cases. In order to prove Theorem B, we will first show $E''(x) = \{0\}$. We begin with some properties of $E'(x)$.

For any $v \in T_xM$, we say that $v$ satisfies condition $(**)$ if

$$(**) \quad \begin{cases} \text{For all vector fields } X_1, \cdots, X_k \text{ defined in a neighborhood} \\ \text{of } x, \text{ each of them lying in a curvature distribution} \\ \text{with curvature normal not equal to zero,} \\ \nabla_v[X_1, [X_2, \cdots [X_{k-1}, X_k] \cdots ]](x) \in \overline{E'(x)}. \end{cases}$$

LEMMA 5.1. *Let $g_t$ be a one-parameter group of isometries of $V$ which preserves $M$. Let $\gamma(t) = g_t(p)$ for some $p \in M$. Then $\overline{E'}|_\gamma$ and $E''|_\gamma$ are smooth vector bundles over $\gamma$ and their direct sum is equal to $TM|_\gamma$. If in addition $\dot\gamma(0)$ satisfies condition $(**)$, then $\nabla_{\dot\gamma(0)}$ preserves this splitting of $TM|_\gamma$.*

*Proof.* For every $t$, since $g_t$ is an extrinsic isometry of $M$, $(g_t)_*$ always maps curvature distributions with nonzero curvature normals to curvature distributions with nonzero curvature normals. Therefore $(g_t)_*(E'(p)) = E'(\gamma(t))$ by the definition of $E'$. Consequently, $(g_t)_*(\overline{E'(p)}) = \overline{E'(\gamma(t))}$ and $(g_t)_*(E''(p)) = E''(\gamma(t))$. This shows that $\overline{E'}|_\gamma$ and $E''|_\gamma$ are smooth vector bundles over $\gamma$. It then follows immediately from the definition that the direct sum of these two bundles is $TM|_\gamma$.

Assume that $\dot\gamma(0)$ satisfies condition $(**)$. To prove the second assertion of the lemma, we only need to prove that $\nabla_{\dot\gamma(0)}$ preserves one of those two bundles, say $E''|_\gamma$. We need to show that if $w$ is a locally defined smooth section of $E''|_\gamma$, then $\nabla_{\dot\gamma(0)}w$ is perpendicular to $\overline{E'(p)}$. Since $E'(p)$ is dense in $\overline{E'(p)}$, it suffices to show that if $u$ is an arbitrary iterated bracket of vectors lying in curvature distributions with nonzero curvature normals, then $u(p) \perp \nabla_{\dot\gamma(0)}w$. In fact $\left\langle \nabla_{\dot\gamma(0)}w, u(p) \right\rangle = -\left\langle w(p), \nabla_{\dot\gamma(0)}u \right\rangle$ which is zero because $\dot\gamma(0)$ satisfies condition $(**)$. □

Actually, condition $(**)$ is satisfied by all tangent vectors of $M$. To prove this, we need the following interesting relationship between the bracket operation and the covariant derivative on $M$.

LEMMA 5.2. *For any vector field $X$ on $M$, denote the orthogonal projection of $X$ to a curvature distribution $E_k$ by $(X)_{E_k}$. If $E_i$ and $E_j$ are two distinct curvature distributions, then for any vector fields $X_i \in E_i$ and $X_j \in E_j$,*

$$(\nabla_{X_i} X_j)_{E_k} = c(i,j,k)([X_i, X_j])_{E_k}$$

*where $c(i,j,k)$ is a constant which can be computed explicitly in terms of curvature normals $n_i$, $n_j$, and $n_k$ by the following formula*

$$c(i,j,k) = \frac{\langle n_i - n_j, n_i - n_k \rangle}{\|n_i - n_j\|^2}.$$

*Proof.* By Lemma 2.1 of [HL], for any vector field $X_k \in E_k$,

$$\langle \nabla_{X_i} X_j, X_k \rangle (n_j - n_k) = \langle \nabla_{X_j} X_i, X_k \rangle (n_i - n_k).$$

Subtracting $\langle \nabla_{X_i} X_j, X_k \rangle (n_i - n_k)$ from both sides, we have

$$\langle \nabla_{X_i} X_j, X_k \rangle (n_j - n_i) = -\langle [X_i, X_j], X_k \rangle (n_i - n_k).$$

Taking the inner product with $-(n_i - n_j)$ and dividing both sides by $\|n_i - n_j\|^2$, we have

$$\langle \nabla_{X_i} X_j, X_k \rangle = \frac{\langle n_i - n_j, n_i - n_k \rangle}{\|n_i - n_j\|^2} \langle [X_i, X_j], X_k \rangle.$$

This proves the lemma. □

LEMMA 5.3. *Every $v \in T_x M$ satisfies condition $(**)$.*

*Proof.* We first observe that, since $\overline{\bigoplus_{n_j \neq 0} E_j(x)} \subset \overline{E'(x)}$, for any $X \in T_x M$, $X \in \overline{E'(x)}$ if and only if $(X)_{E_0} \in \overline{E'(x)}$. The rest of the proof is divided into two steps. Throughout the proof, whenever we say an iterated bracket, we always mean an iterated bracket of some locally defined vector fields, each of which lies in a curvature distribution with nonzero curvature normal. We will also keep in mind Lemma 2.2 of [HL] in our calculations below.

*Step* 1. If $v_i \in E_i(x)$ where the corresponding curvature normal $n_i \neq 0$, then $v_i$ satisfies condition $(**)$.

Extend $v_i$ to a smooth vector field in $E_i$. Let $X_j$ be a vector field lying in $E_j$ where $n_j \neq 0$. If $n_i = n_j$, then $\nabla_{v_i} X_j \in E_i$ since each leaf of $E_i$ is totally geodesic. If $n_j \neq n_i$, then by Lemma 5.2, $(\nabla_{v_i} X_j)_{E_0} = c(i,j,0)([v_i, X_j])_{E_0} \in \overline{E'}$. Therefore for any vector field $X \perp E_0$, $\nabla_{v_i} X \in \overline{E'}$. It also follows from the definition that for such a vector field $X$, $[v_i, X] \in \overline{E'}$.

Let $w$ be an iterated bracket of length $k$. Write $w = w_0 + w_1$ where $w_0 \in E_0$ and $w_1 \perp E_0$. Since $\nabla_{v_i} w_1 \in \overline{E'}$, we only need to show that $\nabla_{v_i} w_0 \in \overline{E'}$. In fact, by Lemma 5.2,

$$(\nabla_{v_i} w_0)_{E_0} = c(i,0,0)([v_i, w_0])_{E_0} = c(i,0,0)([v_i, w] - [v_i, w_1])_{E_0}.$$

This vector lies in $\overline{E'}$ since the first term in the last expression is an iterated bracket of length $k+1$.

*Step* 2. Every $v \in T_x M$ satisfies condition $(**)$.

Extend $v$ to a smooth vector field in a neighborhood of $x$ and write $v = \sum_{j \in I} v_j$ where $I$ is the index set of the set of all curvature normals and $v_j \in E_j$ for all $j \in I$. Let $w$ be an iterated bracket. We will prove $\nabla_v w \in \overline{E'}$ by induction on the length of $w$.

For the base case where the length of $w$ is 0 (i.e. $k=1$), $w$ is a vector field in some curvature distribution $E_i$ with curvature normal $n_i \neq 0$. By Step 1, we only need to show that $\nabla_{v_0} w \in \overline{E'}$. In fact, by Lemma 5.2, $(\nabla_{v_0} w)_{E_0} = 0$ since $c(0,i,0) = 0$. This proves the base case.

Assume the lemma is true for all iterated brackets of length less than or equal to $k$. Let $w$ be an iterated bracket of length $k+1$. We can write $w = [X_i, w']$ where $X_i$ lies in some curvature distribution $E_i$ with $n_i \neq 0$ and $w'$ is an iterated bracket of length $k$. We also write $w' = w_0 + w_1$ where $w_0 \in E_0$ and $w_1 \perp E_0$. Now

$$\begin{aligned}
\nabla_v w &= \nabla_v \nabla_{X_i} w' - \nabla_v \nabla_{w'} X_i \\
&= R_{v,X_i} w' + \nabla_{X_i} \nabla_v w' + \nabla_{[v,X_i]} w' - \nabla_v \nabla_{w_0} X_i - \nabla_v \nabla_{w_1} X_i,
\end{aligned}$$

where $R$ is the curvature operator of $M$. By the Gauss equation, $(R_{Z_1,Z_2} Z_3)_{E_0} = 0$ for all vectors $Z_j \in T_x M$, $j = 1,2,3$. In particular $R_{v,X_i} w' \in \overline{E'}$. By the induction hypothesis, both $\nabla_v w'$ and $\nabla_{[v,X_i]} w'$ lie in $\overline{E'}$. By Corollary 4.2 and Lemma 4.5, $X_i$ is tangent to an orbit of a one-parameter group of extrinsic isometries on $M$. Combining results of Step 1 and Lemma 5.1, we have $\nabla_{X_i} \nabla_v w' \in \overline{E'}$. By Lemma 5.2, $\nabla_{w_0} X_i \perp E_0$ since $c(0,i,0) = 0$. Therefore $\nabla_{w_0} X_i$ is a sum of iterated brackets of length 0. It then follows from the base case that $\nabla_v \nabla_{w_0} X_i \in \overline{E'}$. It remains to show that $\nabla_v \nabla_{w_1} X_i \in \overline{E'}$. In fact

$$\nabla_v \nabla_{w_1} X_i = R_{v,w_1} X_i + \nabla_{w_1} \nabla_v X_i + \nabla_{[v,w_1]} X_i.$$

The first term on the right side lies in $\overline{E'}$ since $R_{v,w_1} X_i \perp E_0$. It follows from the base case that both $\nabla_{[v,w_1]} X_i$ and $\nabla_v X_i$ lie in $\overline{E'}$. Since $w_1$ can be written as a sum of possibly infinitely many vectors each of which lies in a curvature distribution with nonzero curvature normal, results of Step 1 and Lemma 5.1 imply that $\nabla_{w_1} \nabla_v X_i \in \overline{E'}$. This finishes the proof of the lemma. □

LEMMA 5.4. *Let $g_t$ be a one-parameter group of isometries of $V$ which preserve $M$. Let $\gamma(t) = g_t(p)$ for some $p \in M$. If $v$ is a parallel tangent vector field along $\gamma$ such that $v(p) \in E''(p)$, then $v$ is constant along $\gamma$ and $v(p) \in E''(\gamma(t))$ for all $t$.*

*Proof.* Let $v(t) = v(\gamma(t))$. We write $v = v_0 + v_1$ where $v_1(t) \in \overline{E'(\gamma(t))}$ and $v_0(t) \in E''(\gamma(t))$. By Lemma 5.1, $v_0$ and $v_1$ are smooth sections of the bundles $E''|_\gamma$ and $\overline{E'}|_\gamma$ respectively. Moreover, by Lemmas 5.3 and 5.1, $\nabla_{\dot\gamma(t)} v_0(t) \in E''(\gamma(t))$ and $\nabla_{\dot\gamma(t)} v_1(t) \in \overline{E'(\gamma(t))}$. This implies that both $v_0$ and $v_1$ are parallel tangent vector fields along $\gamma$ because $v$ is parallel along $\gamma$. By assumption, $v(p) = v_0(p)$. Therefore $v(t) = v_0(t)$ for all $t$ by the uniqueness of parallel vector fields. This proves that $v(t) \in E''(\gamma(t))$ for all $t$.

It remains to show that $v$ is constant along $\gamma$. In fact, since $v$ is parallel along $\gamma$, $\frac{d}{dt} v(t) \in \nu_{\gamma(t)} M$. On the other hand, for any normal vector $n \in \nu_{\gamma(t)} M$, we have

$$\left\langle \frac{d}{dt} v(t),\ n \right\rangle = \langle A_n v(t),\ \dot\gamma(t) \rangle = 0,$$

where $A$ is the shape operator of $M$. The last equality is due to the fact that $v(t) \in E''(\gamma(t)) \subset E_0(\gamma(t))$. Consequently $\frac{d}{dt} v(t) \equiv 0$. The lemma is thus proved. □

To continue the proof of Theorem B, we need to introduce two new equivalence relations. For any points $x$ and $y$ in $M$, we say $x \approx y$ if there exists a piecewise differentiable curve connecting $x$ and $y$ such that each differentiable piece of this curve is a part of the orbit of a one-parameter group of isometries of $V$ which preserve $M$. We say $x \approx_0 y$ if there exists a sequence of points $x_k \in M$, $k = 0, 1, \cdots, n$, such that $x_0 = x$, $x_n = y$, and for every $k = 1, \cdots, n$, either $x_k \approx x_{k-1}$ or $x_k \in x_{k-1} + E_0(x_{k-1})$. For any $x \in M$, we define

$$\mathfrak{O}(x) := \{y \in M \mid y \approx x\}, \quad \text{and} \quad \mathfrak{O}_0(x) := \{y \in M \mid y \approx_0 x\}.$$

By Corollary 4.2 and Lemma 4.5, we have:

LEMMA 5.5. *For any $p \in M$, $Q(p) \subset \mathfrak{O}(p)$ and $Q_0(p) \subset \mathfrak{O}_0(p)$, where $Q(p)$ is the equivalence class defined in the introduction and $Q_0(p)$ is defined in Section 3.*

It follows from Theorem D that $Q_0(p)$ is dense in $M$ for all $p \in M$. Therefore we have:

COROLLARY 5.6. *For any $p \in M$, $\overline{\mathfrak{O}_0(p)} = M$.*

Moreover, we also have the following relation between $\mathfrak{O}_0(p)$ and $\mathfrak{O}(p)$.

LEMMA 5.7. *For any $p \in M$, $\mathfrak{O}_0(p) = \bigcup_{x \in \mathfrak{O}(p)} (x + E_0(x))$.*

*Remark.* We do not know whether the corresponding relation between $Q_0(p)$ and $Q(p)$ holds. In fact this is the main reason for introducing the classes $\mathfrak{O}_0(p)$ and $\mathfrak{O}(p)$.

*Proof.* Let $g_t$ be a one-parameter group of isometries on $V$ which preserve $M$. Since every extrinsic isometry of $M$ preserves the zero curvature distribution on $M$, for any $x \in M$ and $y \in x + E_0(x)$, $g_t(y) \in g_t(x + E_0(x)) = g_t(x) + E_0(g_t(x))$ for all $t$. By definition, $g_t(x) \in \mathfrak{O}(x)$. The lemma then follows from a trivial induction argument. □

PROPOSITION 5.8. *For any $x \in M$, $E'(x)$ is dense in $T_x M$.*

*Proof.* We only need to prove that $E''(x) = \{0\}$. In fact, if there exists $v \in E''(x)$, then by Lemma 5.4, $v \in E''(y)$ for all $y \in \mathfrak{O}(x)$. Since $E''(y) \subset E_0(y)$, for any $z \in y + E_0(y)$, $v \in E_0(z) = E_0(y)$. By Lemma 5.7 and Corollary 5.6 and the continuity of $E_0$, we have $v \in E_0(z)$ for all $z \in M$. Therefore the straight line $z + \mathbb{R} \cdot v$ is contained in $M$ for all $z \in M$. Consequently, if $v \neq 0$, then $M$ is a product of the straight line $x + \mathbb{R} \cdot v$ and a proper submanifold of $M$. This contradicts the assumption that $M$ is irreducible. Hence we have $E''(x) = \{0\}$. □

Now we are ready to prove the main theorems.

*Proof of Theorem* B. Let $D$ be the set of all locally defined vector fields on $M$ which lie in some curvature distributions with nonzero curvature normals. Then the set of reachable points of $D$ starting from $p$ is just $Q(p)$. By Proposition 5.8, $D^*(x) = E'(x)$ is dense in $T_x M$ for all $x \in M$. Therefore the theorem follows from Theorem D. □

*Proof of Theorem* A. This theorem follows from Proposition 4.4 and Theorem B. □

*Proof of Theorem* C. For any $y \in L_P(x)$, by Theorem A, there exists an extrinsic isometry, say $f$, of $M$ which maps $x$ to $y$. Moreover, by Theorem 4.1 and the proof of Proposition 4.4, $f$ can be chosen to preserve curvature distributions and curvature normals. Therefore $f_*$ preserves the distribution $D_P$. Consequently, $f(L_P(x)) = L_P(y) = L_P(x)$ and $f(W_P(x)) = W_P(y) = W_P(x)$. This finishes the proof of the theorem. □

## Appendix. An infinite dimensional version of Chow's theorem

In this appendix, we give the proof of Theorem D. We begin with a lemma whose proof is essentially due to Jens Heber.

LEMMA A.1. *Let $V$ be a Hilbert space and $B \subset V$ be a closed subset with $B \neq \emptyset, V$. Then there exists a truncated cone*

$$C := C(a, v, \alpha, r) := \{x \in V \mid x \neq a, \ \angle(x - a, v) \leq \alpha, \ \langle x - a, v \rangle \leq r\} \cup \{a\},$$

where $a, v \in V$, $\|v\| = 1$, $\alpha \in (0, \pi/2)$ and $r > 0$, such that $C \cap B = \{a\}$.

*Proof.* Let $a_1 \in B$ and $x \in V \setminus B$. Put $r_1 = \|x - a_1\|$ and $v = (1/r_1)(x - a_1)$. Then $x \in C(a_1, v, \alpha, r_1)$ for any $\alpha > 0$. Since $V \setminus B$ is open, there exists $\alpha \in (0, \pi/2)$ and $\rho_1 \in [0, r_1)$, such that

$$C(a_1, v, \alpha, r_1) \cap B \subset C(a_1, v, \alpha, \rho_1).$$

We will fix $v$ and $\alpha$ in the following and put $C(a, r) := C(a, v, \alpha, r)$. Note that for any $b \in C(a, r)$,

$$C(a, r) \cap C(b, \infty) = C(b, r - \langle b - a, v \rangle).$$

We may assume that $\rho_1$ is chosen to be optimal, i.e. that

$$\rho_1 = \inf\{\rho \geq 0 \mid C(a_1, r_1) \cap B \subset C(a_1, \rho)\}.$$

We define inductively a sequence of nested truncated cones

$$C(a_1, r_1) \supset C(a_2, r_2) \supset \cdots$$

and real numbers $\rho_1, \rho_2, \ldots$ by choosing $a_{k+1} \in C(a_k, \rho_k) \cap B$ with $\langle a_{k+1} - a_k, v \rangle \geq \rho_k/2$ and putting $r_{k+1} = r_k - \langle a_{k+1} - a_k, v \rangle$ and

$$\rho_{k+1} = \inf\{\rho \geq 0 \mid C(a_{k+1}, r_{k+1}) \cap B \subset C(a_{k+1}, \rho)\}.$$

It then follows that

$$C(a_{k+1}, r_{k+1}) = C(a_k, r_k) \cap C(a_{k+1}, \infty) \subset C(a_k, r_k),$$

and $C(a_{k+1}, \rho_{k+1}) \subset C(a_k, \rho_k)$. Therefore $\rho_{k+1} \leq \rho_k - \langle a_{k+1} - a_k, v \rangle \leq \rho_k/2$ and $r_{k+1} - \rho_{k+1} \geq r_k - \rho_k$. Hence $\lim_{k \mapsto \infty} \rho_k = 0$ and $r_k$ converges to some $r > 0$, since $r_{k+1} \geq r_{k+1} - \rho_{k+1} \geq \cdots \geq r_1 - \rho_1 > 0$. Because for any positive integer $l$, $a_{k+l} \in C(a_{k+l-1}, \rho_{k+l-1}) \subset C(a_k, \rho_k)$,

$$\|a_{k+l} - a_k\| \leq \operatorname{diam} C(a_k, \rho_k) = \rho_k \cdot \operatorname{diam} C(a, 1)$$

for any $a \in V$. Therefore the $a_k$, $k = 1, 2, \ldots$, form a Cauchy sequence and thus converge to some point $a \in B$. The lemma now follows from

$$C(a, r) \cap B = \bigcap_{k=1}^{\infty} C(a_k, r_k) \cap B \subset \bigcap_{k=1}^{\infty} C(a_k, \rho_k) = \{a\}. \qquad \square$$

*Remark.* If $\dim V < \infty$, a proof follows trivially by taking a ball in the complement of $B$ of maximal radius and a point of $B$ on the boundary of the ball.

Let $B$ be an arbitrary subset of a Hilbert manifold $N$. For any $x \in B$ and $v \in T_x N$, we say that $v$ is tangent to $B$ if there exists a curve $\gamma$ in $N$ such that $\gamma(0) = x$, $\dot\gamma(0)$ exists and is equal to $v$, and $\gamma(t) \in B$ for all $t$. We denote the set of all tangent vectors to $B$ at $x$ by $T_x B$.

COROLLARY A.2. *Let $N$ be a connected Hilbert manifold and $B \subset N$ be a closed nonempty subset. If for every $x \in B$, $T_xB$ is dense in $T_xN$, then $B = N$.*

*Proof.* We only need to show that $B$ is also open. In fact if $B$ is not open, then there exists a boundary point $a$ of $B$. Let $(\phi, U)$ be a chart around $a$ with $\phi(U) = V$, where $V$ is the Hilbert space on which $N$ is modeled. Then $\phi(B \cap U)$ is a closed proper subset of $V$. By Lemma A.1, there exists a truncated cone in $V$ which meets $\phi(B \cap U)$ precisely in its vertex. This contradicts the assumption that $T_xB$ is dense in $T_xN$ for all $x \in N$. □

*Remark.* The proof actually shows that if $B$ is a proper closed subset of $N$, then the set of $x \in \partial B$ which admit a truncated cone $C$ with $B \cap C = \{x\}$ is dense in $\partial B$, where we understand by a truncated cone in $N$ a subset diffeomorphic to a truncated cone in a Hilbert space.

For the proof of Theorem D, we need another ingredient, which is probably well known, at least in the finite dimensional case. Since we cannot find any reference, we provide a proof here for completeness.

LEMMA A.3. *Let $X_1, \cdots, X_k$ be vector fields on a Hilbert manifold $N$, and $p_0 \in N$. Then there exists a neighborhood $U$ of $p_0$, an $\varepsilon > 0$ and for all $t \in (-\varepsilon, \varepsilon)$, a (local) diffeomorphism $\alpha_t$ which is a composition of the (local) one-parameter groups $\phi_t^i$ of $X_i$ such that for any $p \in U$ and any differentiable function $f$,*

$$\left.\frac{d^l}{dt^l}\right|_{t=0} f(\alpha_t p) = \begin{cases} 0 & \text{if } l = 1, \cdots, k-1, \\ (k!)[X_1, \cdots, [X_{k-1}, X_k] \cdots]_p(f) & \text{if } l = k. \end{cases}$$

*Proof.* Since the statement is a local one, we may assume for simplicity that all $\phi_t^i$, $i = 1, \cdots, k$, are defined on $N$ and for all $t$. If $k = 1$, we let $\alpha_t = \phi_t^1$. Now we proceed by induction on $k$. If $\alpha_t = \alpha_t^k$ has already been constructed for $X_1, \cdots, X_k$ with the above properties and $\phi_t$ is the one-parameter group of a further vector field $X$, we define $\alpha_t^{k+1} = \alpha_t^{-1} \phi_t^{-1} \alpha_t \phi_t$ (generalizing the well known definition for $k = 1$). From

$$0 = \left.\frac{d^l}{dt^l}\right|_{t=0} f(\alpha_t^{-1} \alpha_t p) = \sum_{l_1+l_3=l} \frac{l!}{l_1! l_3!} \left.\frac{\partial^l}{\partial t_1^{l_1} \partial t_3^{l_3}}\right|_{(0,0)} f(\alpha_{t_1}^{-1} \alpha_{t_3} p)$$

for all $l \geq 1$, we have

$$\left.\frac{d^l}{dt^l}\right|_{t=0} f(\alpha_t^{-1} p) = \begin{cases} 0 & \text{if } l = 1, \cdots, k-1, \\ -\left.\frac{d^k}{dt^k}\right|_{t=0} f(\alpha_t p) & \text{if } l = k, \end{cases}$$

for all $p \in N$ and any smooth function $f$ on $N$. Therefore

$$\left.\frac{\partial^l}{\partial t_1^{l_1} \partial t_2^{l_2} \partial t_3^{l_3} \partial t_4^{l_4}}\right|_{(0,0,0,0)} f(\alpha_{t_1}^{-1} \phi_{t_2}^{-1} \alpha_{t_3} \phi_{t_4} p) = 0$$

if $1 \leq l_1 \leq k-1$ or $1 \leq l_3 \leq k-1$. Furthermore, since

$$\frac{d^l}{dt^l}\bigg|_{t=0} f(\alpha_t^{k+1} p)$$
$$= \sum_{l_1+l_2+l_3+l_4=l} \frac{l!}{l_1!l_2!l_3!l_4!} \frac{\partial^l}{\partial t_1^{l_1} \partial t_2^{l_2} \partial t_3^{l_3} \partial t_4^{l_4}}\bigg|_{(0,0,0,0)} f(\alpha_{t_1}^{-1} \phi_{t_2}^{-1} \alpha_{t_3} \phi_{t_4} p)$$

and

$$0 = \frac{d^l}{dt^l}\bigg|_{t=0} f(\phi_t^{-1} \phi_t p) = \sum_{l_2+l_4=l} \frac{l!}{l_2!l_4!} \frac{\partial^l}{\partial t_2^{l_2} \partial t_4^{l_4}}\bigg|_{(0,0)} f(\phi_{t_2}^{-1} \phi_{t_4} p),$$

for all $l \geq 1$, we get

$$\frac{d^l}{dt^l}\bigg|_{t=0} f(\alpha_t^{k+1} p) = 0 \quad \text{if } 1 \leq l \leq k,$$

and

$$\frac{d^{k+1}}{dt^{k+1}}\bigg|_{t=0} f(\alpha_t^{k+1} p)$$
$$= (k+1) \frac{\partial^{k+1}}{\partial s \partial t^k}\bigg|_{(0,0)} \{f(\alpha_t^{-1}\phi_s^{-1} p) + f(\alpha_t^{-1}\phi_s p) + f(\phi_s^{-1}\alpha_t p) + f(\alpha_t \phi_s p)\}$$
$$= (k+1) \frac{\partial^{k+1}}{\partial s \partial t^k}\bigg|_{(0,0)} \{f(\alpha_t \phi_s p) - f(\phi_s \alpha_t p)\}$$
$$= (k+1)![X, [X_1, \cdots, [X_{k-1}, X_k] \cdots]]_p(f). \qquad \square$$

COROLLARY A.4. *Let $D$ be a set of vector fields which are defined on open subsets of $N$. Then $D^*(x) \subset T_x \Omega_D(x)$ for any $x \in N$.*

*Proof.* If $v$ is an iterated bracket of $X_1, \cdots, X_k \in D$ which is defined in an open neighborhood of $x \in N$, by Lemma A.3, there exists a one-parameter family of local diffeomorphisms $\alpha_t$, which is defined in an open neighborhood $U$ of $x$, such that for all $p \in U$ and $|t| < \varepsilon$, $\alpha_t p \in \Omega_D(p)$ and for any smooth function $f$ on $U$,

$$\frac{d^l}{dt^l}\bigg|_{t=0} f(\alpha_t p) = \begin{cases} 0 & \text{if } l = 1, \cdots, k-1, \\ (k!)v(p)f & \text{if } l = k. \end{cases}$$

Let $\beta_t = \alpha_{\sqrt[k]{t}}$ for $t \geq 0$. Then for every $p \in U$ and $t$ sufficiently small, $\beta_t p \in \Omega_D(p)$ and $\frac{d}{dt}\big|_{t=0} (\beta_t p) = v(p)$, as can be seen by taking the Taylor series of $\alpha_t p$ in a local coordinate system.

For any $v \in D^*$ which is defined in an open neighborhood of $x \in N$, we can write $v = v_1 + \cdots + v_l$ where each $v_i$ is an iterated bracket of vector fields lying in $D$. Let $\beta_t^i$ be the one-parameter family of local diffeomorphisms constructed for $v_i$ as above, and $\gamma(t) = \beta_t^1 \cdots \beta_t^l x$ for small $t$. Then $\gamma(t) \in \Omega_D(x)$ for all $t$ where $\gamma$ is defined and $\dot{\gamma}(0) = v(x)$. Therefore $v(x) \in T_x \Omega_D(x)$. $\square$

*Proof of Theorem* D. We first claim that for any $p \in N$ and $x \in \overline{\Omega_D(p)}$, $\Omega_D(x) \subset \overline{\Omega_D(p)}$. In fact, we may assume that $x_k$, $k = 1, 2, \cdots$, is a sequence of points in $\Omega_D(p)$ which converges to $x$. Suppose that $y \in \Omega_D(x)$ can be joined to $x$ by an integral curve $\gamma$ of a vector field $X \in D$ which is defined on an open subset $U$ of $N$. For sufficiently large $k$, $x_k \in U$. Therefore there is an integral curve of $X$ passing through $x_k$, which is denoted by $\gamma_{x_k}$. Since an integral curve of a smooth vector field depends continuously on its initial point, $\gamma_{x_k}$ converges to $\gamma$. Therefore there is a sequence of points $y_k = \gamma_{x_k}(t_k)$, which is defined for large $k$ and real numbers $t_k$, such that $\lim_{k \mapsto \infty} y_k = y$. Since $y_k \in \Omega_D(x_k) = \Omega_D(p)$, we have $y \in \overline{\Omega_D(p)}$. The claim then follows from a trivial induction argument. By Corollary A.4, $D^*(x) \subset T_x \Omega_D(x) \subset T_x \overline{\Omega_D(p)}$ for all $x \in \overline{\Omega_D(p)}$. The theorem then follows from Corollary A.2. □

Universität Augsburg
Universitätsstrasse, 14
D-86159 Augsburg, Germany
*E-mail address*: ernst.heintze@math.uni-augsburg.de

Max-Planck-Institut für Mathematik
Gottfried-Claren-Strasse 26
D-53225 Bonn, Germany
*Current address*:
Department of Mathematics
Massachussetts Institute of Technology
Cambridge, MA 02139
*E-mail address*: xbliu@math.mit.edu